\documentclass{gtpart}
\usepackage{pinlabel}
\usepackage{upgreek}
\usepackage{hyperref}
\hypersetup{backref,colorlinks=true,citecolor=blue,linkcolor=blue}
\usepackage{calc}
\newlength{\wcwidth}
\newlength{\wcheight}
\newcommand{\widecheck}[1]{\ensuremath{
\settowidth{\wcwidth}{#1}
\settoheight{\wcheight}{#1}
\addtolength{\wcheight}{1pt}
\makebox[0cm][l]{%
\raisebox{\depth+\wcheight}[0cm][0cm]{%
\scalebox{-1}{$\widehat{\hphantom{#1}}$}}}#1
\rule{0pt}{\wcheight+2.5pt}}}
\title{Seiberg-Witten Floer homology and\\symplectic forms on $\mathrm{S^1\times M^3}$}
\author{Cagatay Kutluhan}
\givenname{\c{C}a\u{g}atay}
\surname{Kutluhan}
\address{Department of Mathematics\\ The University of Michigan\\\newline 530 Church Street\\ Ann Arbor, MI 48109 USA}
\email{kutluhan@umich.edu}
\author{Clifford Henry Taubes}
\givenname{Clifford Henry}
\surname{Taubes}
\address{Department of Mathematics\\ Harvard University\\\newline One Oxford Street\\ Cambridge, MA 02138 USA}
\email{chtaubes@math.harvard.edu}
\keyword{Seiberg-Witten Floer homology}
\keyword{symplectic form}
\subject{primary}{msc2000}{57R17}
\subject{primary}{msc2000}{57R57}
\arxivreference{math.SG/08041371}
\volumenumber{}
\issuenumber{}
\publicationyear{}
\papernumber{}
\startpage{}
\endpage{}
\doi{}
\MR{}
\Zbl{}
\received{}
\revised{}
\accepted{}
\published{}
\publishedonline{}
\proposed{}
\seconded{}
\corresponding{}
\editor{}
\version{}
\newtheorem{theorem}{Theorem}[section]

\newtheorem{lemma}[theorem]{Lemma}

\newtheorem{prop}[theorem]{Proposition}

\newtheorem*{mt}{Main Theorem}
\newcommand{\QED}{\hfill$\Box$\medskip}
\newcommand{\sint}{\mathrm{I_{i}}}
\newcommand{\bint}{\mathrm{J_{i,i+1}}}
\newcommand{\Aii}{\mathrm{A_{i,i+1}}}
\newcommand{\Aiit}{{\mathrm{A_{i,i+1}}|_t}}
\newcommand{\aii}{\mathrm{a_{i,i+1}}}
\newcommand{\daii}{\mathrm{\dot{a}_{i,i+1}}}
\newcommand{\Ac}{\mathrm{A}_{\underline{\mathbb{C}}}}

\newcommand{\FAc}{\mathrm{F}_{\mathrm{A}_{\underline{\mathbb{C}}}}}
\newcommand{\FAct}{\mathrm{F}_{\mathrm{A}_{\underline{\mathbb{C}}}|_t}}
\newcommand{\dAc}{\dot{\mathrm{A}}_{\underline{\mathbb{C}}}}
\numberwithin{equation}{section}
\providecommand{\bysame}{\leavevmode\hbox to3em{\hrulefill}\thinspace}
\begin{document}
\begin{abstract}
Let $\mathrm{M}$ be a closed, connected, orientable three-manifold. The purpose of this paper is to study the Seiberg-Witten Floer homology of $\mathrm{M}$ given that $\mathrm{S}^1\times \mathrm{M}$ admits a symplectic form.\end{abstract}
\maketitle
\section{Introduction}
\label{S1}
Suppose $\mathrm{M}$ is a closed, connected, orientable three-manifold such that the product four-manifold $\mathrm{S}^1\times \mathrm{M}$ admits a symplectic form. Let $\upomega$ denote a symplectic form on $\mathrm{S}^1\times \mathrm{M}$. Then, one can write $\upomega$ as
\begin{equation}
\label{1.1}
\upomega=\mathrm{dt}\wedge \upnu + \upmu
\end{equation}
where $\mathrm{dt}$ is a nowhere vanishing 1-form on $\mathrm{S^1}$, $\upnu$ is a section over $\mathrm{S}^1\times \mathrm{M}$ of $\mathrm{T^*M}$ and $\upmu$ is a section over $\mathrm{S}^1\times \mathrm{M}$ of $\mathrm{\wedge^2\;T^*M}$. Let $\mathrm{d}$ denote the exterior derivative along $\mathrm{M}$ factor of $\mathrm{S^1\times M}$. 
Since $\upomega$ is a closed 2-form, one has $\frac{\partial}{\partial t} \upmu= \mathrm{d}\upnu$ and $\mathrm{d}\upmu = 0$. Thus, $\upmu$ is a closed form on $\mathrm{M}$ at any given $t\in\mathrm{S}^1$. Its cohomology class in $\mathrm{H}^2(\mathrm{M};\mathbb{R})$ is denoted by $[\upmu]$. As explained momentarily, the class $[\upmu]$ is non-zero. To see why this is the case, first use the K\"{u}nneth formula to write $\mathrm{H^2}(\mathrm{S^1\times M};\mathbb{R})$ as the direct sum $[\mathrm{dt}]\cup\mathrm{H^1}(\mathrm{M};\mathbb{R})\oplus \mathrm{H^2}(\mathrm{M};\mathbb{R})$ where $[\mathrm{dt}]$ denotes the cohomology class of the 1-form $\mathrm{dt}$. Let $[\upomega]$ denote the cohomology class of the symplectic form $\upomega$. This class appears in the K\"{u}nneth decomposition as $[\mathrm{dt}]\cup[\bar{\upnu}]+[\upmu]$ where $\bar{\upnu}$ is the push-forward from $\mathrm{S^1\times M}$ of the 2-form $\mathrm{dt}\wedge\upnu$. This understood, neither $[\bar{\upnu}]$ nor $[\upmu]$ are zero by virtue of the fact that $[\upomega]\cup[\upomega]$ is non-zero.
\par Our convention is to orient $\mathrm{S}^1$ by $\mathrm{dt}$, and $\mathrm{S^1\times M}$ by $\upomega\wedge\upomega$. Doing so finds that $\upnu\wedge\upmu$ is nowhere zero and so orients $\mathrm{M}$ at any given $t\in\mathrm{S}^1$.
\par Now, fix a $t$-independent Riemannian metric, $\mathfrak{g}$, on $\mathrm{M}$, and let $\ast$ denote the corresponding Hodge star operator. At each $t\in\mathrm{S}^1$ , the 1-form $\ast\upmu$ is a nowhere vanishing 1-form on $\mathrm{M}$ and so defines a homotopy class of oriented 2-plane fields by its kernel.  This 2-plane field is denoted in what follows by $\mathrm{K^{-1}}$. This bundle is oriented by $\upmu$ and so has a corresponding Euler class which we write as $-c_1(\mathrm{K})\in\mathrm{H}^2(\mathrm{M};\mathbb{Z})$.
\par Fix a $\mathrm{spin^c}$ structure on $\mathrm{M}$ and let $\mathbb{S}$ denote the associated spinor bundle, this a Hermitian $\mathbb{C}^2$-bundle over $\mathrm{M}$. At any $t\in\mathrm{S}^1$, the eigenbundles for Clifford multiplication by $\ast\upmu$ on $\mathbb{S}$ split $\mathbb{S}$ as a direct sum, $\mathbb{S}=\mathrm{E}\oplus \mathrm{EK^{-1}}$, where $\mathrm{E}$ is a complex line bundle over $\mathrm{M}$. Here, our convention is to write the $+i|\upmu|$ eigenbundle on the left. The \emph{canonical $spin^c$} structure is that with $\mathrm{E}=\underline{\mathbb{C}}$, the trivial complex line bundle. We use $det(\mathbb{S})$ to denote the complex line bundle $\wedge ^2\mathbb{S}=\mathrm{E^2K^{-1}}$ over $\mathrm{M}$. Note that the assignment of $c_1(\mathrm{E})\in\mathrm{H}^2(\mathrm{M};\mathbb{Z})$ to a given $\mathrm{spin^c}$ structure identifies the set of equivalence classes of $\mathrm{spin^c}$ structures over $\mathrm{M}$ with $\mathrm{H}^2(\mathrm{M};\mathbb{Z})$. This classification of the $\mathrm{spin^c}$ structures over $\mathrm{M}$ is independent of the choice of $t\in\mathrm{S}^1$. For any given class $e\in\mathrm{H^2}(\mathrm{M};\mathbb{Z})$, we use $\mathfrak{s}_e$ to denote the corresponding $\mathrm{spin^c}$ structure. Thus the spinor bundle $\mathbb{S}$ for $\mathfrak{s}_e$ splits as $\mathrm{E}\oplus\mathrm{EK^{-1}}$ with $c_1(\mathrm{E})=e$.
\par P. B. Kronheimer and T. S. Mrowka in \cite{km} associate three versions of the Seiberg-Witten Floer homology to any given $\mathrm{spin^c}$ structure. With $e\in\mathrm{H^2}(\mathrm{M};\mathbb{Z})$ given, the three versions of the Seiberg-Witten Floer homology for the $\mathrm{spin^c}$ structure $\mathfrak{s}_e$ are denoted by Kronheimer and Mrowka and in what follows by $\overline{HM}(\mathrm{M},\mathfrak{s}_e)$, $\widehat{HM}(\mathrm{M},\mathfrak{s}_e)$ and $\widecheck{HM}(\mathrm{M},\mathfrak{s}_e)$. Each of these is a $\mathbb{Z}/p\mathbb{Z}$ graded module over $\mathbb{Z}$ with $p$ the greatest divisor in $\mathrm{H^2}(\mathrm{M};\mathbb{Z})$ of the cohomology class $2e-c_1(\mathrm{K})$, which is the first Chern class of the corresponding version of $\mathbb{S}$. Each of these modules is a $\mathrm{C}^{\infty}$ invariant of $\mathrm{M}$.
\par The purpose of this paper is to prove the following theorem.
\begin{mt}
Let $\mathrm{M}$ be a closed, connected, orientable three-manifold. Suppose that $\mathrm{S^1\times M}$ has the symplectic form $\upomega=\mathrm{dt}\wedge \upnu + \upmu$.  Fix a class $e\in\mathrm{H}^2(\mathrm{M};\mathbb{Z})$ with $2e-c_1(\mathrm{K})=\uplambda[\upmu]$ in $\mathrm{H^2}(\mathrm{M};\mathbb{R})$ for some $\uplambda<0$. Let $\mathfrak{s}_e$ denote the $\mathrm{spin^c}$ structure  corresponding to $e$ via the correspondence defined above. Then  $\overline{HM}(\mathrm{M},\mathfrak{s}_e)$ vanishes, $\widehat{HM}(\mathrm{M},\mathfrak{s}_e)\cong\widecheck{HM}(\mathrm{M},\mathfrak{s}_e)$, and the following hold:
\begin{itemize}
\item If $e=0$, then $\widecheck{HM}(\mathrm{M},\mathfrak{s}_e)\cong\mathbb{Z}$.
\item Suppose $e\neq 0$. Then $\widecheck{HM}(\mathrm{M},\mathfrak{s}_e)$ vanishes if the pull-back of $e$ by the obvious projection map from $\mathrm{S^1\times M}$ onto $\mathrm{M}$ has non-positive pairing with the Poincar\'e dual of $[\upomega]$. \end{itemize} 
\end{mt}
\noindent We say that the \emph{monotonicity condition} is satisfied by a given $\mathrm{spin^c}$ structure $\mathfrak{s}_e$ when $2e-c_1(\mathrm{K})=\uplambda[\upmu]$ holds in $\mathrm{H^2}(\mathrm{M};\mathbb{R})$ for some $\uplambda<0$.
\par As it turns out, our Main Theorem also describes Seiberg-Witten Floer homology for $\mathrm{spin^c}$ structures with $2e-c_1(\mathrm{K})=\uplambda[\upmu]$ in $\mathrm{H^2}(\mathrm{M};\mathbb{R})$ for some $\uplambda>0$. Here is why: Let $e\in\mathrm{H^2}(\mathrm{M};\mathbb{Z})$ be given. Then Proposition 25.5.5 in \cite{km} describes an isomorphism between Seiberg-Witten Floer homology groups for $\mathfrak{s}_e$ and those for $\mathfrak{s}_{c_1(\mathrm{K})-e}$. In particular, if $2e-c_1(\mathrm{K})=\uplambda[\upmu]$ with $\uplambda>0$, then the monotonicity condition is satisfied for the $\mathrm{spin^c}$ structure $\mathfrak{s}_{c_1(\mathrm{K})-e}$ and our Main Theorem applies.
\par The following remarks are meant to give some context to this theorem.  
First, the Euler characteristic of the Seiberg-Witten Floer homology for any given $\mathrm{spin^c}$ structure is called the Seiberg-Witten invariant of the $\mathrm{spin^c}$ structure. Our Main Theorem is consistent with what \cite{t0} claims about Seiberg-Witten invariants of $\mathrm{M}$.
\par Second, suppose that $\mathrm{M}$ fibers over the circle.  Let $f: \mathrm{M}\rightarrow \mathrm{S^1}$ 
denote a locally trivial fibration.  Then, $\mathrm{M}$ admits a metric that makes $f$ harmonic. In this case, the 
pull-back, $df$, by $f$ of the Euclidean 1-form on $\mathrm{S^1} = \mathbb{R}/2\pi \mathbb{Z}$ is a harmonic 1-form.  Hence, the 2-form $\upomega = \mathrm{dt}\wedge df+\ast df$ is symplectic on $\mathrm{S^1\times M}$. When the fiber of $f$ has genus 2 or greater, the monotonicity condition for any $e\in\mathrm{H}^2(\mathrm{M};\mathbb{Z})$ with $e=\upkappa[\ast df]$ for some $\upkappa\leq 0$ is satisfied and the conclusions of our Main Theorem are known to be true. 
\par The third remark concerns the following question: If $\mathrm{S^1\times M}$ admits a symplectic 
form, does $\mathrm{M}$ fiber over $\mathrm{S^1}$? A very recent preprint by S. Friedl and S. Vidussi \cite{fv} asserts an affirmative answer to this qestion. Our Main Theorem with Theorem 1 of Y. Ni in \cite{n} (see also \cite{km2}) gives a different proof that $\mathrm{M}$ fibers over $\mathrm{S^1}$ in the case when $\mathrm{M}$ has first Betti number 1 and $c_1(\mathrm{K})$ is not torsion.
\begin{theorem}
\label{1.1}
Let $\mathrm{M}$ be a closed, connected, irreducible, orientable, three-manifold with first Betti number equal to 1. Let $\upomega$ denote a symplectic form on $\mathrm{S^1\times M}$ such that $c_1(\mathrm{K})$ is not torsion. Then $\mathrm{M}$ fibers over $\mathrm{S^1}$.
\end{theorem}
Note that if $c_1(\mathrm{K})$ is not torsion in $\mathrm{H^2}(\mathrm{M};\mathbb{Z})$, then $c_1(\mathrm{K})=\uplambda[\upmu]$ in $\mathrm{H}^2(\mathrm{M};\mathbb{R})$ with $\uplambda>0$. To see why, let $\upkappa$ denote the cup product pairing between $c_1(\mathrm{K})$ and $[\upomega]$. This has the same sign as $\uplambda$. If $\upkappa<0$, then it follows from \cite{l} or \cite{oo} that $\mathrm{M=S^1\times S^2}$. On the other hand, if $c_1(\mathrm{K})$ is torsion, then it follows from our Main Theorem, Proposition 25.5.5 and Theorem 41.5.2 in \cite{km} that $\mathrm{M}$ has vanishing Thurston (semi)-norm. It follows from a theorem of J. D. McCarthy \cite{mc} with G. Perelman's proof of the Geometrization 
Conjecture that $\mathrm{S^1\times M}$ has a symplectic form in the case when $\mathrm{M}$ is reducible if and 
only if $\mathrm{M}=\mathrm{S^1\times S^2}$. 
\\
\par\noindent\textbf{Proof of Theorem \ref{1.1}.} Let $\mathbf{S}$ denote the generator of $\mathrm{H_2}(\mathrm{M};\mathbb{Z})$ with the property that $\langle c_1(\mathrm{K}),\mathbf{S}\rangle>0$. Note that such a class exists by virtue of the fact noted above that $c_1(\mathrm{K})=\uplambda[\upmu]$ with $\uplambda>0$. Let $\Sigma$ denote a closed, connected, oriented and genus minimizing representative for the class $\mathbf{S}$. Use $g$ to denote the genus of $\Sigma$. It is a consequence of Corollary 40.1.2 in \cite{km} (the adjunction inequality) that $2g-2\geq\langle c_1(\mathrm{K}),\mathbf{S}\rangle$. This is to say that $c_1(\mathrm{K})$ lies in the unit ball as defined by the dual of the Thurston (semi)-norm on $\mathrm{H^2}(\mathrm{M};\mathbb{Z})/Tor$. In fact, $c_1(\mathrm{K})$ is an extremal point in this ball, which is to say that $\langle c_1(\mathrm{K}),\mathbf{S}\rangle=2g-2$. Here is why: our Main Theorem in the present context says that 
$$\bigoplus_{e\in\mathrm{H^2}(\mathrm{M};\mathbb{Z})\;:\; \langle e,\mathbf{S}\rangle<0}\widecheck{HM}(\mathrm{M},\mathfrak{s}_e)\cong\{0\},$$
$$\bigoplus_{e\in\mathrm{H^2}(\mathrm{M};\mathbb{Z})\;:\; \langle e,\mathbf{S}\rangle=0}\widecheck{HM}(\mathrm{M},\mathfrak{s}_e)\cong\mathbb{Z}.$$
Meanwhile, Proposition 25.5.5 in \cite{km} asserts isomorphisms between the Seiberg-Witten Floer homology groups for the $\mathrm{spin^c}$ structure $\mathfrak{s}_e$ and those for the $\mathrm{spin^c}$ structure $\mathfrak{s}_{c_1(\mathrm{K})-e}$. Thus, our Main Theorem also finds that 
\begin{eqnarray}
\label{e1.1}
\nonumber\bigoplus_{e\in\mathrm{H^2}(\mathrm{M};\mathbb{Z})\;:\; \langle e,\mathbf{S}\rangle>\langle c_1(\mathrm{K}),\mathbf{S}\rangle}\widecheck{HM}(\mathrm{M},\mathfrak{s}_e)\cong\{0\},\\
\bigoplus_{e\in\mathrm{H^2}(\mathrm{M};\mathbb{Z})\;:\; \langle e,\mathbf{S}\rangle=\langle c_1(\mathrm{K}),\mathbf{S}\rangle}\widecheck{HM}(\mathrm{M},\mathfrak{s}_e)\cong\mathbb{Z}.
\end{eqnarray}
These last results with Theorem 41.5.2 in \cite{km} imply that $c_1(\mathrm{K})$ is an extremal point of the unit ball as defined by the dual of the Thurston (semi)-norm, that is to say $\langle c_1(\mathrm{K}),\mathbf{S}\rangle=2g-2$. Given (\ref{e1.1}), the assertion made by Theorem \ref{1.1} follows directly from Theorem 1 in \cite{n}.\QED

\par\noindent\textbf{Acknowledgments.} The first author would like to thank his thesis advisor Prof. Daniel Burns for his support throughout the course of this project. He would also like to thank University of Michigan Mathematics Department for their support during the term of Winter 2007. The first author dedicates this result to his parents. The second author is supported in part by the National Science Foundation.
\section{Background on Seiberg-Witten Theory}
\label{S2}
In this section, we present a brief introduction to the theory of Seiberg-Witten invariants of three-manifolds and the monopole Floer homology as defined in the book by Kronheimer and Mrowka \cite{km}. In what follows, $\mathrm{M}$ is a given closed, oriented three-manifold.
\subsection{Algebraic preliminaries}\hspace{1in}\par There is a unique connected double cover of the group $\mathrm{SO}(3)$, namely the group $\mathrm{Spin}(3)=\mathrm{SU}(2)$. The group $\mathrm{Spin^c}(3)$ is defined as the quotient of $\mathrm{U(1)}\times \mathrm{Spin}(3)$ by the diagonal action of $\mathbb{Z}_2$, thus the group $\mathrm{U(2)}$. Fix a Riemannian metric on $\mathrm{M}$. A $\mathrm{spin^c}$ structure on $\mathrm{M}$ can be viewed as a principal $\mathrm{U}(2)$-bundle $\mathrm{\tilde{P}}$ such that $\mathrm{\tilde{P}}\times_{\rho}\mathrm{SO}(3)\cong\mathrm{P}_{\mathrm{SO}(3)}$, the principal $\mathrm{SO}(3)$-bundle associated to the tangent bundle of $\mathrm{M}$.  Here, $\rho$ denotes the natural projection of $\mathrm{U}(2)$ onto $\mathrm{U}(2)/\mathrm{U}(1)=\mathrm{SO}(3)$.
\par A $\mathrm{spin^c}$ structure on $\mathrm{M}$ has an associated Hermitian $\mathbb{C}^2$-bundle, this defined by the defining representation of $\mathrm{U}(2)$. This bundle is denoted by $\mathbb{S}$ and it is called the spinor bundle. Its sections are called spinors. There exists the Clifford algebra homomorphism $\mathfrak{cl}:\wedge \mathrm{T_{\mathbb{C}}^*M}\rightarrow \mathrm{End_{\mathbb{C}}(\mathbb{S})}$ that gives a representation of the bundle of Clifford algebras. 
\par There is also a map $det:\mathrm{U}(2)\rightarrow \mathrm{U(1)}$ defined by the determinant. This representation of $\mathrm{U}(2)$ yields a principal $\mathrm{U(1)}$-bundle $\mathrm{\tilde{P}}\times_{det}\mathrm{U(1)}$. The complex line bundle associated to $\mathrm{\tilde{P}}\times_{det}\mathrm{U(1)}$ is called the determinant bundle of the $\mathrm{spin^c}$ structure, which we denote by $det(\mathbb{S})$, because this line bundle is the second exterior power of the bundle $\mathbb{S}$.
\par The existence of $\mathrm{spin^c}$ structures on $\mathrm{M}$ follows immediately from the fact that $\mathrm{M}$ is parallelizable.  The set of $\mathrm{spin^c}$ structures on $\mathrm{M}$ form a principle bundle over a point for the additive group $\mathrm{H^2}(\mathrm{M};\mathbb{Z})$.  To elaborate, a given cohomology class acts on a given $\mathrm{spin^c}$ structure in such a way that the spinor bundle for the new $\mathrm{spin^c}$ structure is obtained from that of the original one by tensoring with a complex line bundle whose first Chern class is the given class in $\mathrm{H^2(M;\mathbb{Z})}$.
\subsection{Seiberg-Witten Floer homology}\hspace{1in}
\par Let $\mathcal{S}$ denote the set of $\mathrm{spin^c}$ structures on $\mathrm{M}$. A unitary connection $\mathbb{A}$ on $det(\mathbb{S})$ together with the Levi-Civita connection on the orthonormal frame bundle of M determines a $\mathrm{spin^c}$ connection $\mathbf{A}$ on the spinor bundle $\mathbb{S}$. Then the Seiberg-Witten monopole equations are
\begin{eqnarray}\label{2.1}
\nonumber\ast\mathrm{F}_{\mathbb{A}}&=&\uppsi^{\dagger} \uptau \uppsi-i\varrho\\
\mathcal{D}_{\mathbb{A}}\uppsi&=&0.
\end{eqnarray}
Here, the notation is as follows: First, $\mathrm{F}_{\mathbb{A}}\in \Omega^2(\mathrm{M},i \mathbb{R})$ denotes the curvature of the connection $\mathbb{A}$. Second, $\uppsi $ is a section of the spinor bundle $\mathbb{S}$. Third, $\uppsi^{\dagger}\uptau\uppsi$ denotes the section of $i\mathrm{T^{\ast}M}$ which is the metric dual of the homomorphism $\uppsi^{\dagger}\mathfrak{cl(\cdot)}\uppsi:\mathrm{T^{\ast}M}\rightarrow i\mathbb{R}$. Fourth, $\mathcal{D}_{\mathbb{A}}$ is the Dirac operator associated to $\mathbf{A}$, which is defined by
$$\Gamma(\mathbb{S})\stackrel{\nabla_{\mathbf{A}}}\longrightarrow \Gamma(\mathrm{T^*M\otimes \mathbb{S}})\stackrel{\mathfrak{cl}}{\longrightarrow}\Gamma(\mathbb{S}).$$
Finally, $\varrho$ is a fixed smooth co-closed 1-from on $\mathrm{M}$.
\par The equations in (\ref{2.1}) are the variational equations of a functional defined on the configuration space $\mathcal{C}=\mathrm{Conn}(det(\mathbb{S}))\times \mathrm{C}^{\infty}(\mathrm{M};\mathbb{S})$ as
\begin{equation}
\nonumber\mathfrak{csd}(\mathbb{A},\uppsi)=-\frac{1}{2}\int_{\mathrm{M}}\mathbb{(A-A_S)}\wedge (\mathrm{F}_{\mathbb{A}}+\mathrm{F}_{\mathbb{A_S}})-i\int_{\mathrm{M}}\mathbb{(A-A_S)}\wedge\ast\varrho+\int_{\mathrm{M}}\uppsi^{\dagger}\mathcal{D}_{\mathbb{A}}\uppsi.
\end{equation}
Here, $\mathbb{A_S}$ is any given connection fixed in advance on $det(\mathbb{S})$. This is the so-called \emph{Chern-Simons-Dirac} functional.
\par The group of gauge transformations of a $\mathrm{spin^c}$ structure, namely the \emph{gauge group} $\mathcal{G}=\mathrm{C}^{\infty}(\mathrm{M},\mathrm{S^1})$, acts on the configuration space as
\begin{eqnarray}
\nonumber\mathcal{G}\times\mathcal{C}&\longrightarrow& \mathcal{C}\\
\nonumber(u,(\mathbb{A},\uppsi))&\longmapsto&(\mathbb{A}-2u^{-1}\mathrm{d}u, u\uppsi).
\end{eqnarray}
\par The equations in (\ref{2.1}) are invariant under the action of the gauge group. Therefore, one can define the space of equivalence classes of solutions of these equations under the action of the gauge group. This is called the \emph{moduli space}, which we denote by $\mathcal{M}$. The solutions of the equations in (\ref{2.1}) which are of the form $(\mathbb{A},0)$ are called reducible solutions because the stabilizer under the action of the gauge group is not trivial. Solutions with non-zero spinor component are called irreducible. We let $\mathcal{B=C/G}$. It is possible to prove that $\mathcal{M}$ is a sequentially compact subset of $\mathcal{B}$. The gauge group $\mathcal{G}$ acts freely on the space of irreducible solutions of the equations in 
(\ref{2.1}).  If $\varrho$ is suitably generic, then the quotient of this space by $\mathcal{G}$ is a finite set of points in $\mathcal{B}$.
\par To elaborate, let $\underline{\mathbb{R}}$ denote the trivial line bundle over $\mathrm{M}$. Each $(\mathbb{A},\uppsi)\in\mathcal{C}$ has an associated linear operator $\mathcal{L}_{(\mathbb{A},\uppsi)}$ that maps $\mathrm{C^{\infty}}(\mathrm{M};i\mathrm{T^*M\oplus\mathbb{S}}\oplus i\underline{\mathbb{R}})$ onto itself. It is defined as
\begin{equation}
\nonumber\mathcal{L}_{(\mathbb{A},\uppsi)}(\mathrm{b},\upphi,g)=\left(\begin{array}{c}\ast\mathrm{db}-\mathrm{d}g-(\uppsi^{\dagger}\uptau\upphi+\upphi^{\dagger}\uptau\uppsi)\\ \mathcal{D}_{\mathbb{A}}\upphi+\frac{1}{2}\mathfrak{cl}(\mathrm{b})\uppsi+g\uppsi\\ -\mathrm{d}^{\ast}\mathrm{b}-\frac{1}{2}(\upphi^{\dagger}\uppsi-\uppsi^{\dagger}\upphi)\end{array}\right).
\end{equation}
This operator extends to $\mathrm{L^{2}}(\mathrm{M};i\mathrm{T^*M\oplus\mathbb{S}}\oplus i\underline{\mathbb{R}})$ as an unbounded, self-adjoint Fredholm operator with dense domain $\mathrm{{L^{2}}_1}(\mathrm{M};i\mathrm{T^*M\oplus\mathbb{S}}\oplus i\underline{\mathbb{R}})$. It has a discrete spectrum that is unbounded from above and below. The spectrum has no accumulation points, and each eigenvalue has finite multiplicity. 
\par An irreducible solution of the equations in (\ref{2.1}) is called non-degenerate if the kernel of $\mathcal{L}$ is trivial. A generic choice for $\varrho$ renders all such solutions non-degenerate.  In this case, irreducible solutions of the equations in (\ref{2.1}) define isolated points in $\mathcal{B}$. 
\par Seiberg-Witten Floer homology is an infinite dimensional version of the Morse homology theory where $\mathcal{B}$ plays the role of the ambient manifold and the Chern-Simons-Dirac functional plays the role of the ``Morse'' function. As the critical points of the Chern-Simons-Dirac functional are solutions of the equations in (\ref{2.1}), the latter are used, as in Morse theory, to label generators of the chain complex. The analog of a non-degenerate critical 
point is a solution of the equations in (\ref{2.1}) whose version of $\mathcal{L}$ has trivial kernel. Here, the point is that $\mathcal{L}$ is, formally, the Hessian of the Chern-Simons-Dirac functional. 
\par As the Hessian in finite dimensional Morse theory can be used to define the 
grading of the Morse complex, it is also the case here that the operator $\mathcal{L}$ is used to define 
a grading for each generator of the Seiberg-Witten Floer homology chain complex.  In 
particular, $\mathcal{L}$ can be used to associate an integer degree to each non-degenerate solution of the equations in
(\ref{2.1}), in fact, to any given pair in $\mathcal{C}$ whose version of $\mathcal{L}$ has trivial kernel. It is enough to say here that this degree involves the notion of spectral 
flow for families of self adjoint operators such as $\mathcal{L}$.  In general, only the $mod(p)$ reduction of this degree is gauge invariant, where $p$ is the greatest integer divisor of 
$c_1(det(\mathbb{S}))$.   
\par The analog in this context of a gradient flow line in finite dimensional Morse 
theory is a smooth map $s\mapsto (\mathbb{A}(s),\uppsi(s))$ from $\mathbb{R}$ into $\mathcal{C}$ that obeys the rule
\begin{eqnarray}
\nonumber \frac{\partial}{\partial s}\mathbb{A} &=& -\ast\mathrm{F_{\mathbb{A}}}+\uppsi^{\dagger} \uptau \uppsi-i\varrho\\
\nonumber\frac{\partial}{\partial s}\uppsi &=& -\mathcal{D}_{\mathbb{A}}\uppsi.
\end{eqnarray}
This can also be written as $\frac{\partial}{\partial s}(\mathbb{A},\uppsi)=-\nabla_{\mathrm{L}^2}\mathfrak{csd}|_{(\mathbb{A},\uppsi)}$ where $\nabla_{\mathrm{L}^2}$ denotes the $\mathrm{L}^2$-gradient of $\mathfrak{csd}$. An \emph{instanton} is a solution of these equations on $\mathbb{R}\times\mathrm{M}$ that converges to a solution of the equations in (\ref{2.1}) on each end as $|s|$ tends to infinity.
\par The differential on the Seiberg-Witten Floer homology chain complex is defined using a suitably perturbed version of these instanton equations.  As in finite dimensional Morse theory, a perturbation is in general necessary in order to have a well defined count of solutions.  The perturbed equations can be viewed as defining the analog of what in finite dimensions would be the equations that define the flow lines of a pseudo-gradient vector field for the given function.  Kronheimer and Mrowka describe in Chapter III of their book \cite{km} a suitable Banach space, $\mathcal{P}$, of such perturbations.  Kronheimer and Mrowka prove that there is a residual set of such perturbations with the following properties:  Each can be viewed as perturbations of $\mathfrak{csd}$, in which case the resulting version of (\ref{2.1}) can serve to define generators of the Seiberg-Witten Floer homology chain complex.  Meanwhile, the resulting instanton equations can serve to define the differential on this chain complex.
\par Note for future reference that $\mathcal{P}$ contains a subspace, $\Omega$, of 1-forms $\varrho$ for use in (\ref{2.1}).  The induced norm on $\Omega$ dominates all of the $\mathrm{C^k}$-norms on $\mathrm{C^{\infty}(M; T^*M)}$.  In fact, if M is assumed to have a real analytic structure, then each $\varrho\in\Omega$ is itself real analytic. An important point to note later on is that the function $\mathfrak{csd}$ decreases along any solution of its gradient flow equations.  This is also the case for the just described perturbed analog of $\mathfrak{csd}$ and the solutions of the latter's gradient flow equations.
\par  \section{Outline of the Proof}
\label{S3}
\par Our purpose in this section is to outline our proof of the Main Theorem and in doing so, state the principle analytic results we will need.  The proofs for most of the assertions made in this section are deferred to the subsequent sections of this article.
\par Fix $t\in\mathrm{S^1}$, and let $\mathrm{M}_t$ denote the slice $\mathrm{M}_t = \{t\}\times\mathrm{M}$. A version of the Seiberg-Witten equations on $\mathrm{M}_t$ can be defined as follows: Let $\upvarpi_{\mathbb{S}}$ be the harmonic 2-form on $\mathrm{M}$ representing the class $2\pi c_1(det(\mathbb{S}))$.  Fix a connection, $\mathbb{A_S}$, on $det(\mathbb{S})$ with curvature 2-form $-i\upvarpi_{{\mathbb{S}}}$. Then, any given connection on $det(\mathbb{S})$ is of the form $\mathbb{A_S}+2\mathrm{a}$ for $\mathrm{a}\in \mathrm{C}^{\infty}(\mathrm{M};i\mathrm{T^*M})$. 
\par Now, fix $r\geq 1$ and $t\in\mathrm{S}^1$. We consider the equations
\begin{eqnarray}
\label{3.1}
\nonumber\ast\mathrm{da}&=&r(\uppsi^{\dagger}\uptau\uppsi-i\ast\upmu)+\frac{i}{2}\ast\upvarpi_{\mathbb{S}}\\
\mathcal{D}_{\mathbb{A}}\uppsi&=&0,
\end{eqnarray}
where $\upmu$ is the 2-form defined by the symplectic form. Suitably rescaling $\uppsi$, we see that these are a version of the equations in (\ref{2.1}). These equations are the variational equations of a functional defined as
\begin{equation}
\label{3.2}\mathfrak{a}(\mathbb{A_S}+2\mathrm{a},\uppsi)=-\frac{1}{2}\int_{\mathrm{M}_t}\mathrm{a}\wedge (\mathrm{da}-i\upvarpi_{\mathbb{S}})-ir\int_{\mathrm{M}_t}\mathrm{a}\wedge \upmu+r\int_{\mathrm{M}_t}\uppsi^{\dagger}\mathcal{D}_{\mathbb{A}}\uppsi,\end{equation}
where $\mathrm{a}\in\mathrm{C}^{\infty}(\mathrm{M};i\mathrm{T^{\ast}M})$ and $\uppsi\in\mathrm{C}^{\infty}(\mathrm{M};\mathbb{S})$.
\par For future purposes, we introduce a new functional on $\mathcal{C}$. Fix $r\geq1$, $t\in\mathrm{S}^1$ and for $(\mathbb{A},\uppsi)\in \mathcal{C}$ let
\begin{equation}
\label{3.31}
\mathcal{E}(\mathbb{A},\uppsi)=i\int_{\mathrm{M}_t}\upnu\wedge\mathrm{da}.
\end{equation}
\par Our approach is to consider $\mathrm{S^1}\times\mathrm{M}$ as a 1-parameter family of three-dimensional manifolds, each a copy of $\mathrm{M}$ and parametrized by $t\in\mathrm{S^1}$.  We use the gauge equivalence classes of solutions of the equations in (\ref{3.1}) on $\mathrm{M}_t$ (when non-degenerate) to define the generators of the Seiberg-Witten Floer homology. Here it is important to remark that the solutions of the equations in (\ref{3.1}) can serve this purpose for any $r\geq 1$ because we assume that $c_1(det(\mathbb{S}))=\uplambda[\upmu]$ with $\uplambda<0$. For the same reason, (\ref{3.1}) has no reducible solutions.
\par Here, we remark that what is written in (\ref{3.1}) has  \emph{period class} $-[\upmu]$ in the sense of \cite{km}. The assumption that $[\upmu]$ is a negative multiple of $c_1(det(\mathbb{S}))$ is what is called the \emph{monotone} case in \cite{km}. As is explained in Chapter VIII of \cite{km}, the results from the case of \emph{exact} perturbations carry onto the monotone case almost without any change, and there are canonical isomorphisms between the Floer homology groups defined here and the relevant Seiberg-Witten Floer homology groups.
\par There is one more important point to make here: The only $t$-dependence in (\ref{3.1}) is due to the appearance of the 2-form $\upmu$ through the latter's $t$-dependence on $t\in\mathrm{S}^1$.
to define generators of the corresponding Seiberg-Witten Floer homology.  Note that the $t$-dependence is due entirely to the appearance of the 2-form $\upmu$ and its dependence on $t$. 
\par We suppose our main theorem is false, and hence that there are at least two generators of the Seiberg-Witten Floer homology for each $t\in\mathrm{S^1}$.  Note in this regard that there is at least one generator for the $\mathrm{E} = \underline{\mathbb{C}}$ case because the fact that $\mathrm{S^1\times M}$ is symplectic implies, via the main theorem in \cite{t0}, that the Seiberg-Witten invariant for the canonical $\mathrm{spin^c}$ structure on $\mathrm{S^1\times M}$ is equal to 1.  If there are at least two generators, then there are at least two solutions.  Our plan is to use the large $r$ behavior of at least one of these solutions to construct nonsense from the assumed existence of two or more generators.  
\par What follows describes what we would like to do.  Given the existence of two or more non-zero Seiberg-Witten Floer homology classes, we would like to use a variant of the strategy from \cite{t1} and \cite{t2} to find, for large enough $r\geq 1$ and for each $t\in\mathrm{S}^1$, a set $\Theta_t\subset\mathrm{M}_t$ of the following sort:  $\Theta_t$ is a finite set of pairs of the form $(\upgamma, \mathrm{m})$ with $\upgamma\subset \mathrm{M}_t$ a closed integral curve of the vector field that generates the kernel of $\upmu|_t$, and $\mathrm{m}$ is a positive integer. These are constrained so that no two pair have the same integral curve.  In addition, with each $\upgamma$ oriented by $\ast\upmu|_t$, the formal sum $\Sigma_{(\upgamma,\mathrm{m})\in\Theta_t} \mathrm{m}\upgamma$ represents the Poincar\'e dual to $c_1(\mathrm{E})$ in $\mathrm{H}_1(\mathrm{M}_t; \mathbb{Z})$. We would also like the graph $t\rightarrow\Theta_t$ to sweep out a smooth, oriented surface $\mathrm{S}\subset\mathrm{S^1\times M}$ whose fundamental class gives the Poincar\'e dual to $c_1(\mathrm{E})$ in $\mathrm{H}_2(\mathrm{S^1\times M}; \mathbb{Z})$.  Note in this regard that such a surface is oriented by the vector field $\frac{\partial}{\partial t}$ and by the 1-form $\upnu$ that appears when we write $\upomega= \mathrm{dt}\wedge\upnu+\upmu$.  In particular, $\upomega|_{\mathrm{TS}}$ is positive and so the integral of $\upomega$ over $\mathrm{S}$ is positive.  On the other hand, the integral of $\upomega$ over $\mathrm{S}$ must be non-positive if the cup product of $[\upomega]$ with $c_1(\mathrm{E})$ is non-positive.  This is the fundamental contradiction.
\par As it turns out, we cannot guaranteed that $\Theta_t$ exists for all $t\in\mathrm{S^1}$, only for most $t$, where `most' has a precise measure-theoretic definition.  Even so, we have control over enough of $\mathrm{S^1}$ to obtain a contradiction which is in the spirit of the one described from any violation to the assertion of our main theorem.
\par To elaborate, consider first the existence of $\Theta_t$.  What follows is the key to this existence question.
\begin{prop}\label{p3.1}
Fix a bound on the $\mathrm{C^3}$-norm of $\upmu$, and fix constants $\mathcal{K} > 1$ and $\updelta > 0$.  Then, there exists $\upkappa > 1$ with the following significance:  Suppose that $r \geq\upkappa$, $t\in\mathrm{S}^1$ and $(\mathbb{A},\uppsi)$ is a solution of the $t$ and $r$ version of the equations in (\ref{3.1}) such that $\mathcal{E}(\mathbb{A},\uppsi)\leq \mathcal{K}$ and such that $\mathrm{sup}_{\mathrm{M}} (|\upmu|-|\uppsi| ^2) >\updelta$.  Then there exists a set $\Theta_t$ of the sort described above. 
\end{prop}
\par The next proposition says something about when we can guarantee Proposition \ref{p3.1}'s condition on $|\uppsi|$:
\begin{prop}\label{p3.2}
Fix a bound on the $\mathrm{C}^3$-norm of $\upmu$. Then, there exists $\upkappa>1$ such that if $r\geq\upkappa$, then the following are true:
\begin{itemize}
\item Suppose that $\mathbb{S}=\underline{\mathbb{C}}\oplus\mathrm{K}^{-1}$. Then, for any $t\in\mathrm{S^1}$, there exists a unique gauge equivalence class of solutions $(\mathbb{A}_{\underline{\mathbb{C}}},\uppsi_{\underline{\mathbb{C}}})$ of the $t$ and $r$ version of the equations in (\ref{3.1}) with $|\uppsi_{\underline{\mathbb{C}}}|\geq |\upmu|^{1/2} -\upkappa^{-1}$. Moreover, these solutions are non-degenerate with $|\uppsi_{\underline{\mathbb{C}}}|\geq |\upmu|^{1/2}-\upkappa r^{-1/2}$ and $\mathcal{E}(\mathbb{A}_{\underline{\mathbb{C}}},\uppsi_{\underline{\mathbb{C}}})\leq\upkappa$.
\item Suppose that $\mathbb{S}=\mathrm{E}\oplus\mathrm{EK^{-1}}$ with $c_1(\mathrm{E})\neq 0$. If $(\mathbb{A},\uppsi)$ is a solution of any given $t\in\mathrm{S}^1$ version of the equations in (\ref{3.1}), then there exists points in $\mathrm{M}$ where $|\uppsi|\leq\upkappa r^{-1/2}$.
\end{itemize}
\end{prop}
\par Proposition 3.1 raises the following, perhaps obvious, question:  

\begin{center}\emph{How do we find, other than by Proposition \ref{p3.2}, solutions with $\mathcal{E}$ bounded at large $r$?}\end{center}

\par To say something about this absolutely crucial question, remark that Proposition \ref{p3.1} here has an almost verbatim analog that played a central role in \cite{t1} and \cite{t2}.  These papers use the analog of (\ref{3.1}) with $\ast\upmu$ replaced by a contact 1-form to prove the existence of Reeb vector fields.  The contact 1-form version of $\mathcal{E}$ replaces the form $\upnu$ with the contact 1-form also.  The existence of an $r$-independent bound on the contact 1-form version of $\mathcal{E}$ played a key role in the arguments given in \cite{t1} and \cite{t2}.  The existence of the desired bound on the contact 1-form version of $\mathcal{E}$ exploits the \textbf{$r$-dependence} of the functional $\mathfrak{a}$.

\par We obtain the desired $r$-independent bound on our version of $\mathcal{E}$ for most $t\in\mathrm{S}^1$ by exploiting the \textbf{$t$-dependence} of $\mathfrak{a}$.  To say more about this, it proves useful now to introduce a spectral flow function, $\mathcal{F}$, for certain configurations in $\mathcal{C}$.  There are three parts to its definition.  
Here is the first part:  Fix a section $\uppsi_{\mathrm{E}}$ of $\mathbb{S}$ so that the $(\mathbb{A}_{\mathbb{S}}, \uppsi_{\mathrm{E}})$ version of the operator $\mathcal{L}$ 
as defined in Section \ref{S2} is non-degenerate.  Use $\mathcal{L}_{\mathrm{E}}$ to denote the latter operator.  The second part introduces the version of $\mathcal{L}$ that is relevant to (\ref{3.1}); it is obtained from the 
original by taking into account the rescaling of $\uppsi$.  In particular, it is defined by  
\begin{equation}\label{4.18}
\mathcal{L}_{(\mathbb{A},\uppsi)}(\mathrm{b},\upphi,g)=\left(\begin{array}{c}\ast\mathrm{db}-\mathrm{d}g-2^{-1/2}r^{1/2}(\uppsi^{\dagger}\uptau\upphi+\upphi^{\dagger}\uptau\uppsi)\\ \mathcal{D}_{\mathbb{A}}\upphi+2^{1/2}r^{1/2}(\mathfrak{cl}(\mathrm{b})\uppsi+g\uppsi)\\ -\mathrm{d}^{\ast}\mathrm{b}-2^{-1/2}r^{1/2}(\upphi^{\dagger}\uppsi-\uppsi^{\dagger}\upphi)\end{array}\right)
\end{equation}
for each $(\mathrm{b},\upphi,g)\in\mathrm{C^{\infty}}(\mathrm{M};i\mathrm{T^*M\oplus\mathbb{S}}\oplus i\underline{\mathbb{R}})$. Thus, $\mathcal{L}_{\mathrm{E}}$ is the $r = 1$ version of (\ref{4.18}) as defined using $(\mathbb{A}_{\mathbb{S}}, \uppsi_{\mathrm{E}})$.  To start the third part of 
the definition, suppose that $(\mathbb{A}, \uppsi)\in\mathcal{C}$ is non-degenerate in the sense that the operator 
$\mathcal{L}_{(\mathbb{A}, \uppsi)}$ as depicted in (\ref{4.18}) has trivial kernel.  As explained in \cite{t1} and \cite{t2}, there is a well 
defined spectral flow from the operator $\mathcal{L}_{\mathrm{E}}$ to $\mathcal{L}_{(\mathbb{A}, \uppsi)}$ (see, also \cite{t3}).  This integer is 
the value of $\mathcal{F}$ at $(\mathbb{A}, \uppsi)$.  Note that $\mathcal{F}(\cdot)$ is defined on the complement of a codimension-1 
subvariety in $\mathcal{C}$.  As such, it is piecewise constant. In general, only the $mod(p)$ reduction of $\mathcal{F}$ is gauge invariant where $p$ is the greatest divisor of the class $c_1(det(\mathbb{S}))$.
\par The function $\mathfrak{a}$ is not invariant under the action of $\mathcal{G}$ on $\mathcal{C}$; and, as just noted, neither is $\mathcal{F}$ when 
$c_1(det(\mathbb{S}))$ is non-torsion.  However, our assumption that $c_1(det(\mathbb{S}))=\uplambda[\upmu]$ in $\mathrm{H}^2(\mathrm{M};\mathbb{R})$ implies the 
following:  There exists a constant $\mathfrak{C}$ independent of $r\geq 1$ and $t \in \mathrm{S^1}$ such that  
 $$\mathfrak{a}^{\mathcal{F}}=\mathfrak{a}+r\mathfrak{C}\mathcal{F}$$
is invariant under the action of $\mathcal{G}$. 
To say more about the role of $\mathfrak{a}^{\mathcal{F}}$ requires a digression for two preliminary propositions. They are used to associate a value of $\mathfrak{a}^{\mathcal{F}}$ to each generator of the Seiberg-Witten Floer homology.
\begin{prop}\label{p3.3}
Fix $r\geq 1$ and $\updelta>0$. Then there exist a $t$-independent 1-form $\upsigma\in\Omega$ with $\mathcal{P}$ norm bounded by $\updelta$ such that the following is true: Replace $\upmu$ by $\upmu+\mathrm{d}\upsigma$.
\begin{itemize}
\item The resulting 2-form $\upomega=\mathrm{dt}\wedge\upnu+\upmu$ is symplectic.
\item There exists finite sets $\mathfrak{T}_r$ and ${\mathfrak{T}_r}'$ in $\mathrm{S}^1$ such that if $t\in\mathrm{S}^1\setminus\mathfrak{T}_r$, then $\mathfrak{a}^{\mathcal{F}}$ distinguishes distinct gauge equivalence classes of solutions of the $t$ and $r$ version of the equations in (\ref{3.1}). On the other hand, if $t\in\mathrm{S}^1\setminus{\mathfrak{T}_r}'$ all solutions of the $t$ and $r$ version of the equations in (\ref{3.1}) are non-degenerate.
\item There exists a countable set $\mathfrak{S}_r\in\mathrm{S}^1$ that contains $\mathfrak{T}_r\cup{\mathfrak{T}_r}'$ with accumulation points on the latter such that if $t\in\mathrm{S}^1\setminus\mathfrak{S}_r$, then the gauge equivalence classes of solutions of the equations in (\ref{3.1}) can be used to label the generators of the Seiberg-Witten Floer complex. In this regard, the degree of any generator can be taken to be mod(p) reduction of the negative of the spectral flow function $\mathcal{F}$. 	
\end{itemize}
\end{prop}
\par\noindent\textbf{Proof.} The claim in the first bullet of the proposition is obvious. As for the second and third bullets, the proof of these two follow directly from the arguments used in Sections 2a and 2b 
of \cite{t2}.  The latter prove the analog of the second and third bullets of Proposition \ref{p3.3} 
where $r$ varies rather than $t$. With only notational changes, they also prove the second and 
third bullets here. 
\QED
\par Suppose now that $t\in\mathrm{S}^1\setminus\mathfrak{S}_r$ and that $\uptheta$ is a non-zero Seiberg-Witten Floer homology class. Let $\mathfrak{n}=\Sigma\mathrm{z_\mathfrak{i}}\mathfrak{c_i}$ denote a cycle that represents $\uptheta$ as defined using the $t$ and $r$ version of the equations in (\ref{3.1}). Here $\mathrm{z}_{\mathfrak{i}}\in\mathbb{Z}$ and $\mathfrak{c_i}\in \mathcal{C/G}$ is a gauge equivalence class of solutions of the $t$ and $r$ version of the equations in (\ref{3.1}). Let $\mathfrak{a}^{\mathcal{F}}[\mathfrak{n};t]$ denote the maximum value of $\mathfrak{a}^{\mathcal{F}}$ on the set of generators $\{\mathfrak{c_i}\}$ with $\mathrm{z}_{\mathfrak{i}}\neq0$. Set ${\mathfrak{a}^{\mathcal{F}}}_{\uptheta}$ to denote the minimal value in the resulting set $\{\mathfrak{a}^{\mathcal{F}}[\mathfrak{n};t]\}$.
\begin{prop}\label{p3.4}
The various $t\in\mathrm{S^1}\setminus\mathfrak{S}_r$ versions of the Seiberg-Witten Floer homology groups can be identified in a degree preserving manner so that if $\uptheta $ is any given non-zero class, then the function ${\mathfrak{a}^{\mathcal{F}}}_{\uptheta}(\cdot)$ on $\mathrm{S}^1\setminus\mathfrak{S}_r$ extends to the whole of $\mathrm{S}^1$ as a continuous, Lipschitz function that is smooth on the complement of $\mathfrak{T}_r$. Moreover, if $\mathrm{I}\subset\mathrm{S}^1\setminus\mathfrak{T}_r$ is a component, then there exists $\mathrm{I}'\subset\mathrm{S}^1$ containing the closure of $\mathrm{I}$ and a smooth map $\mathfrak{c}_{\uptheta,\mathrm{I}}:\mathrm{I}'\rightarrow\mathcal{C}$ that solves the corresponding version of the equations in (\ref{3.1}) at each $t\in\mathrm{I}'$ and is such that ${\mathfrak{a}^{\mathcal{F}}}_{\uptheta}(t)=\mathfrak{a}^{\mathcal{F}}(\mathfrak{c}_{\uptheta,\mathrm{I}}(t))$ at each $t\in\mathrm{I}'$.
\end{prop}
\par\noindent\textbf{Proof.} The proof is, but for notational changes and two additional remarks, identical to 
that of Proposition 2.5 in \cite{t2}.  To set the stage for the first remark, fix a base point $0\in\mathrm{S^1}\setminus\mathfrak{S}_r$. The identifications of the Seiberg-Witten Floer homology groups given by 
adapting what is done in \cite{t2} may result in the following situation:  As $t$ increases from $0$, 
these identifications results at $t = 2\pi$ in an automorphism, $\mathrm{U}$, on the $t = 0$ version of the 
Seiberg-Witten Floer homology.  This automorphism need not obey ${\mathfrak{a}^{\mathcal{F}}}_{\mathrm{U}\uptheta}={\mathfrak{a}^{\mathcal{F}}}_{\uptheta}$. If not, 
then it follows using Proposition \ref{p3.3} that the identifications made at $t<2\pi$ to define $\mathrm{U}$ 
can be changed if necessary as $t$ crosses points in $\mathfrak{T}_r$ so that the new version of $\mathrm{U}$ does 
obey ${\mathfrak{a}^{\mathcal{F}}}_{\mathrm{U}\uptheta}={\mathfrak{a}^{\mathcal{F}}}_{\uptheta}$. The second remark concerns the fact that any given $\mathfrak{c}_{\uptheta,\mathrm{I}}$ is unique up to 
gauge equivalence.  This follows from Proposition \ref{p3.3}'s assertion that the function $\mathfrak{a}^{\mathcal{F}}$ distinguishes the Seiberg-Witten solutions when $t\in\mathrm{S}^1\setminus\mathfrak{T}_r$.\QED
\par When $\mathrm{E}=\underline{\mathbb{C}}$, we need to augment what is said in Proposition \ref{p3.4} with the following:
\begin{prop}\label{p3.10}
Suppose that $\mathrm{E}=\underline{\mathbb{C}}$ and that there are at least two non-zero Seibeg-Witten Floer homology classes. Then, the identifications made by Proposition \ref{p3.4} between the various $t\in\mathrm{S}^1$ versions of the Seiberg-Witten Floer homology groups can be assumed to have the following property. There is a non-zero class $\uptheta$ such that none of Proposition \ref{p3.4}'s maps $\mathfrak{c}_{\uptheta,\mathrm{I}}$ send the corresponding interval $\mathrm{I}'$ to a solution in the gauge equivalence class of Proposition \ref{p3.2}'s solution $(\mathbb{A}_{\underline{\mathbb{C}}},\uppsi_{\underline{\mathbb{C}}})$.
\end{prop}
\par\noindent\textbf{Proof.} At any given $t\in\mathrm{S}^1$, there is a class $\uptheta$ with $\mathfrak{c}_{\uptheta,\mathrm{I}}$ not gauge equivalent to $(\mathbb{A}_{\underline{\mathbb{C}}},\uppsi_{\underline{\mathbb{C}}})$. It then follows from Proposition \ref{p3.3} that such is the case for any $t\in\mathrm{S}^1\setminus\mathfrak{T}_r$. This understood, Proposition \ref{p3.4}'s isomorphisms can be changed as $t$ crosses a point in $\mathfrak{T}_r$ while increasing from $t=0$ to insure that no version of $\mathfrak{c}_{\uptheta,\mathrm{I}}$ gives the same gauge equivalence class as $(\mathbb{A}_{\underline{\mathbb{C}}},\uppsi_{\underline{\mathbb{C}}})$. \QED
\par Let $\mathrm{I}$ denote a component of $\mathrm{S}^1\setminus\mathfrak{T}_r$. The assignment of $t\in\mathrm{I}'$ to $\mathcal{E}(\mathfrak{c}_{\uptheta,\mathrm{I}}(\cdot))$ associates to $\uptheta$ a smooth function on $\mathrm{I}'$. View this function on $\mathrm{I}$ as the restriction from $\mathrm{S}^1\setminus\mathfrak{T}_r$ of a function, $\mathcal{E}_{\uptheta}$.  Note that the latter need not extend to $\mathrm{S}^1$ as a continuous function.
\par With the function ${\mathfrak{a}^{\mathcal{F}}}_{\uptheta}$ understood, we come to the heart of the matter, which is the formula for the derivative for this function on any given interval $\mathrm{I}\subset\mathrm{S}^1\setminus\mathfrak{T}_r$: Let $\mathfrak{c}_{\uptheta, \mathrm{I}}$ be as described in Proposition \ref{p3.4}. Then 
\begin{equation}
\label{4.3}
\frac{\mathrm{d}}{\mathrm{d}t}\mathfrak{a}^{\mathcal{F}}(\mathfrak{c}_{\uptheta,\mathrm{I}}(t))=-ir\int_{\mathrm{M}_t}\upnu\wedge\mathrm{da}=-r\mathcal{E}_{\uptheta}.
\end{equation}
To explain, keep in mind that $\mathfrak{c}_{\mathrm{I}}$ is a critical point of $\mathfrak{a}^{\mathcal{F}}$ and so the chain rule for the derivative of $\mathfrak{a}^{\mathcal{F}}(\mathfrak{c}_{\uptheta,\mathrm{I}}(\cdot))$ yields 
\begin{equation}
\label{4.2}
\frac{\mathrm{d}}{\mathrm{d}t}\mathfrak{a}^{\mathcal{F}}(\mathfrak{c}_{\uptheta,\mathrm{I}}(t))=-ir\int_{\mathrm{M}_t}\mathrm{a}\wedge\frac{\partial}{\partial t}\upmu;
\end{equation}
and this is the same as (\ref{4.3}) because $\upomega$ is a closed form. Indeed, write $\upomega=\mathrm{dt}\wedge\upnu+\upmu$ to see that the equation $\mathrm{d}\upomega=0$ requires $\frac{\partial}{\partial t}\upmu=\mathrm{d}\upnu$.  This understood, an integration by parts equates (\ref{4.2}) to (\ref{4.3}).    
\par We get bounds on $\mathcal{E}_{\uptheta}$ after integrating (\ref{4.3}) around $\mathrm{S}^1$. Given that ${\mathfrak{a}^{\mathcal{F}}}_{\uptheta}$  is continuous, integration of the left-hand side over $\mathrm{S^1}$ gives zero. Thus, we conclude that
\begin{equation}\label{3.32}
\int_{\mathrm{S}^1}\mathcal{E}_{\uptheta}=0.
\end{equation}
This formula tells us that $\mathcal{E}_{\uptheta}$ is bounded at some points in $\mathrm{S}^1$. To say more, we use the fact that $\upomega\wedge\upomega>0$ to prove
\begin{lemma}\label{p3.5}
There exists a constant $\upkappa>1$ with the following significance: Suppose that $r\geq\upkappa$, $t\in\mathrm{S}^1$, and $(\mathbb{A},\uppsi)$ is a solution of the corresponding version of the equations in (\ref{3.1}). Then, $\mathcal{E}(\mathbb{A},\uppsi)\geq-\upkappa$.
\end{lemma}
Granted this lower bound on $\mathcal{E}$, the next result follows as a corollary:
\begin{lemma}\label{p3.6}
There exists a constant $\upkappa>1$ with the following significance: Fix $r\geq\upkappa$ so as to define the set $\mathfrak{S}_r\subset\mathrm{S}^1$. Let $\uptheta$ denote a non-zero Seiberg-Witten Floer homology class. Let $n$ denote a positive integer.Then, the measure of the set in $\mathrm{S^1}\setminus\mathfrak{S}_r$ where $\mathcal{E}_{\uptheta}\geq 2^n$ is less than $\upkappa2^{-n}$. 
\end{lemma}
\par\noindent\textbf{Proof.} Given the lower bound provided by Lemma \ref{p3.5}, this follows easily from (\ref{3.32}). \QED

\par Given what has been said so far, we have the desired sets $\Theta_t\subset\mathrm{M}_t$ for points $t$ in the complement of a closed set with non-empty interior in $\mathrm{S^1}$. On the face of it, this is far from what we need, which is a surface $\mathrm{S}\subset \mathrm{S^1\times M}$ that is swept out by such points. As we show below, we can make due with what we have. In particular, we first change our point of view and interpret integration of $\upomega$ over a surface in $\mathrm{S^1\times M}$ as integration over $\mathrm{S^1\times M}$ of the product of $\upomega$ and a closed 2-form $\Phi$ that represents the Poincar\'e dual of the surface. We then construct a 2-form $\Phi$ on $\mathrm{S^1\times M}$ that is localized near the surface swept out by $\uptheta_t$ on most of $\mathrm{S^1\times M}$. This partial localization is enough to prove that $\int_{\mathrm{S^1\times M}}\upomega\wedge\Phi>0$ when this integral should be zero or negative. The existence of such a form gives the nonsense that proves the Main Theorem.
\par The construction of $\Phi$ requires first some elaboration on what is said in Proposition \ref{p3.1}. To set the stage, suppose that $(\mathbb{A},\uppsi)$ is a solution of some $t\in\mathrm{S}^1$ version of the equations in (\ref{3.1}). We will write the section $\uppsi$ of $\mathrm{\mathbb{S}=E\oplus EK^{-1}}$ with respect to the splitting defined by $\ast\upmu|_t$ as $\uppsi=(\upalpha,\upbeta)$ where $\upalpha$ is a section of $\mathrm{E}$ and $\upbeta$ is a section of $\mathrm{EK^{-1}}$. 
\begin{prop}\label{p3.7}
Fix a bound on the $\mathrm{C}^3$-norm of $\upmu$, and fix constants $\mathcal{K}>1$ and $\updelta>0$. There exists $\upkappa>1$ with the following significance: Suppose that $r\geq\upkappa$, $t\in\mathrm{S}^1$, and $(\mathbb{A}=\mathrm{A}_0+2\mathrm{A},\uppsi=(\upalpha,\upbeta))$ is a solution of the equations in (\ref{3.1}) with $\mathcal{E}(\mathbb{A},\uppsi)\leq\mathcal{K}$ and with $\mathrm{sup_M}(|\upmu|-|\uppsi|^2)>\updelta$. Then,
\begin{itemize}
\item There exists a finite set $\Theta_t$  whose typical element is a pair $(\upgamma,\mathrm{m})$ with $\upgamma\subset\mathrm{M}_t$  a closed integral curve tangent to the kernel of $\upmu$, and with $\mathrm{m}$ a positive integer.  Distinct pairs in $\Theta_t$ have distinct curves, and $\Sigma_{(\upgamma,\mathrm{m})\in\Theta_t} \mathrm{m}\upgamma$ generates the Poincar\'e dual to $c_1(\mathrm{E})$ in $\mathrm{H_1}(\mathrm{M}_t;\mathbb{Z})$.
\item Each point where $|\upalpha|^2<|\upmu|-\updelta$ has distance $\upkappa r^{-1/2}$ or less from a curve in $\Theta_t$, and also from some point in $\upalpha^{-1}(0)$.
\item Fix $(\upgamma,\mathrm{m})\in\Theta_t$. Let $\mathrm{D}\subset\mathbb{C}$ denote the closed unit disk centered at the origin and $\upvarphi:\mathrm{D}\rightarrow\mathrm{M}_t$ denote a smooth embedding such that all the points in $\upvarphi(\partial\mathrm{D})$ have distance $\upkappa r^{-1/2}$ or more from any loop in $\Theta_t$. Assume in addition that $\upvarphi(\mathrm{D})$ has intersection 1 with $\upgamma$. Fix a trivialization of the bundle $\upvarphi^{\ast}\mathrm{E}$ over $\mathrm{D}$ so as to view $\upvarphi^{\ast}\upalpha$ as a smooth map from $\mathrm{D}$ into $\mathbb{C}$. The resulting map is non-zero on $\partial\mathrm{D}$ and has degree $\mathrm{m}$ as a map from $\partial\mathrm{D}$ into $\mathbb{C}\setminus\{0\}$.
\end{itemize}
\end{prop}

We now fix $r$ very large so as to define the set $\mathfrak{T}_r=\{t_{\mathrm{i}}\}_{\mathrm{i}=1,..,\mathrm{N}_r}$. We set $t_{\mathrm{N}_r+1}=t_1$ and take the index $i$  to increase in accordance with the orientation of $\mathrm{S}^1$. For each $\mathrm{i}$, we use Propositions \ref{p3.4} and \ref{p3.10} to provide $\mathfrak{c}_{\uptheta,[t_{\mathrm{i}},t_{\mathrm{i+1}}]}$ which we write as $(\mathbb{A}_{\mathrm{i,i+1}},\uppsi_{\mathrm{i,i+1}})$. We view the connection $\mathbb{A}_{\mathrm{i,i+1}}$ as defining a connection on the line bundle $det(\mathbb{S})$ over $\mathrm{I}'\times\mathrm{M}$ where $\mathrm{I}'\in\mathrm{S}^1$ is some open neighborhood of $[t_{\mathrm{i}},t_{\mathrm{i+1}}]$. We also view the $t\in[t_{\mathrm{i}},t_{\mathrm{i+1}}]$ versions of Proposition \ref{p3.2}'s connection $\mathbb{A}_{\underline{\mathbb{C}}}$ as a connection on the bundle $\mathrm{K^{-1}}$ over $[t_{\mathrm{i}},t_{\mathrm{i+1}}]\times\mathrm{M}$. Note in this regard that $\mathrm{K}^{-1}$ is the determinant line bundle for the canonical $\mathrm{spin^c}$ structure with spinor bundle $\mathbb{S}_0=\underline{\mathbb{C}}\oplus\mathrm{K}^{-1}$.

\par With $r$ large and $\updelta >0$ very small, we define $\Phi$ on $[t_{\mathrm{i}}+\updelta,t_{\mathrm{i+1}}-\updelta]\times\mathrm{M}$ to be $\frac{i}{2\pi}(\mathrm{F}_{\mathbb{A}_{\mathrm{i,i+1}}}-\mathrm{F_{\mathbb{A}_{\underline{\mathbb{C}}}}})$. This done, we have yet the task of describing $\Phi$ on the part of $\mathrm{S^1\times M}$ where $t\in[t_{\mathrm{i}}-\updelta,t_{\mathrm{i}}+\updelta]$. We do this as follows: If $\updelta>0$ is sufficiently small, then Proposition \ref{p3.7} asserts that $\mathfrak{c}_{\uptheta,[t_{\mathrm{i}},t_{\mathrm{i+1}}]}$ is defined on the interval $[t_{\mathrm{i}}-\updelta,t_{\mathrm{i+1}}+\updelta]$, and likewise $\mathfrak{c}_{\uptheta,[t_{\mathrm{i-1}},t_{\mathrm{i}}]}$ is defined on the interval $[t_{\mathrm{i-1}}-\updelta,t_{\mathrm{i}}+\updelta]$. This understood, we find a suitable gauge transformations so as to write $\mathbb{A}_{\mathrm{i-1,i}}=\mathbb{A_S}+2\mathrm{a_{i-1,i}}$ and $\mathbb{A}_{\mathrm{i,i+1}}=\mathbb{A_S}+2\mathrm{a_{i,i+1}}$ on $[t_{\mathrm{i}}-\updelta,t_{\mathrm{i}}+\updelta]\times\mathrm{M}$. In particular, these gauge transformations are chosen so that the spectral flow between the respective $(\mathbb{A}_{\mathrm{i-1,i}},\uppsi_{\mathrm{i-1,i}})$ and $(\mathbb{A}_{\mathrm{i,i+1}},\uppsi_{\mathrm{i,i+1}})$ versions of (\ref{4.18}) is zero. We then interpolate between $\mathrm{a}_{\mathrm{i-1,i}}$ and $\mathrm{a}_{\mathrm{i,i+1}}$ on $[t_{\mathrm{i}}-\updelta,t_{\mathrm{i}}+\updelta]\times\mathrm{M}$ using a smooth bump function, $\mathrm{v}$ so as to define a connection $\mathbb{A}_i=\mathbb{A_S}+2(1-\mathrm{v})\mathrm{a}_{\mathrm{i-1,i}}+2\mathrm{v}\mathrm{a}_{\mathrm{i,i+1}}$ on $det(\mathbb{S})$ over $[t_{\mathrm{i}}-\updelta,t_{\mathrm{i}}+\updelta]\times\mathrm{M}$.  With this connection in hand, we define $\Phi$ to be $\frac{i}{2\pi}(\mathrm{F}_{\mathbb{A}_{\mathrm{i}}}-\mathrm{F}_{\mathbb{A}_{\underline{\mathbb{C}}}})$ on $[t_{\mathrm{i}}-\updelta,t_{\mathrm{i}}+\updelta]\times\mathrm{M}$. The continuity of the function $t\rightarrow{\mathfrak{a}^{\mathcal{F}}}_{\uptheta}(t)$ is then used to prove the following:
\begin{prop}\label{p3.9}
Fix a bound on the $\mathrm{C}^3$-norm of $\upmu$. There exists $\upkappa>1$ such that if $r\geq\upkappa$ and if $\updelta>0$ is sufficiently small, then
\begin{itemize}
\item $\Phi$ is twice the first Chern class of a bundle of the form $\mathrm{E\otimes L}$ where $c_1(\mathrm{L})$ has zero cup product with $[\upomega]$.
\item $\int_{\mathrm{S^1}\times\mathrm{M}}\upomega\wedge\Phi>0$.
\end{itemize}
\end{prop}
What is claimed by Proposition \ref{p3.9} is not possible given that the first chern class 
of $\mathrm{E}$ is assumed to have non-positive cup product with the class defined by $\upomega$.  Thus there 
can be no counter example to the claim made by our Main Theorem. 
\section{Analytic Estimates}
\label{S4}
This section contains proofs of Propositions \ref{p3.1} and \ref{p3.2} as well as the proof of Lemma \ref{p3.5}. 
\label{S4}
\par Many of the following arguments in this section exploit two fundamental a priori bounds for solutions of the large $r$ versions of (\ref{3.1}). To start with, write a section $\uppsi$ of $\mathrm{\mathbb{S}=E\oplus EK^{-1}}$ as $\uppsi=(\upalpha,\upbeta)$ where $\upalpha$ is a section of $\mathrm{E}$ and $\upbeta$ is a section of $\mathrm{EK^{-1}}$. Then, the next lemma supplies the fundamental estimates on the norms of $\upalpha$ and $\upbeta$.
\begin{lemma}\label{l4.1} Fix a bound on the $\mathrm{C}^3$-norm of $\upmu$. Then, there are constants $c,c'>0$ with the following significance: Suppose that $(\mathbb{A},\uppsi=(\upalpha,\upbeta))$ is a solution of a given $t\in\mathrm{S}^1$ and $r\geq1$ version of the equations in (\ref{3.1}). Then, 
\begin{itemize}
\item $|\upalpha|\leq |\upmu|^{1/2}+c\;r^{-1}$
\item $|\upbeta|^2\leq c'\;r^{-1}(|\upmu| - |\upalpha|^2)+c\;r^{-2}.$
\end{itemize}
\end{lemma}
\vspace{0.05in}
\par\noindent\textbf{Proof.} This lemma is the same as Lemma 2.2 in \cite{t1} except for the inevitable appearance of $|\upmu|$. We will give the proof in this new context. 
\par Since $\mathcal{D_{\mathbb{A}}}\uppsi=0$, one has $\mathcal{D_{\mathbb{A}}}^2\uppsi=0$ as well. Then, the Weitzenb\"{o}ck formula for $\mathcal{D_{\mathbb{A}}}^2$ yields
\begin{equation}
\label{3.6}
\mathcal{D_{\mathbb{A}}}^2\uppsi=\nabla^{\dagger}\nabla\uppsi+\frac{1}{4}\mathcal{R}\;\uppsi-\frac{1}{2}\mathfrak{cl}(\ast\mathrm{F}_{\mathbb{A}})\uppsi=0
\end{equation}
where $\mathcal{R}$ denotes the scalar curvature of the Riemannian metric. Contract this equation with $\uppsi$ to see that
\begin{equation}
\label{3.7}
\frac{1}{2}\mathrm{d}^{\ast}\mathrm{d}|\uppsi|^2+|\nabla\uppsi|^2+\frac{r}{2}|\uppsi|^2(|\uppsi|^2-|\upmu|-\frac{c_0}{r})\leq 0.
\end{equation}
where $c_0>0$ is a constant depending only on the supremum of $|\upvarpi_{\mathbb{S}}|$ and the infimum of the scalar curvature.
\par Now, introduce $\uppsi=|\upmu|^{1/2}\;\uppsi'$, therefore $\upalpha=|\upmu|^{1/2}\;\upalpha'$ and $\upbeta=|\upmu|^{1/2}\;\upbeta'$. Then, one can rewrite (\ref{3.7}) as follows:
\begin{eqnarray}
\label{3.8}
\nonumber&&\frac{|\upmu|}{2}\mathrm{d}^{\ast}\mathrm{d}|{\uppsi}'|^2-< \mathrm{d}|\upmu|, \mathrm{d} |{\uppsi}'|^2>+\frac{1}{2}|\uppsi'|^2\mathrm{d^{\ast}d|\upmu|}\\
&&+\frac{r}{2}|\upmu\|{\uppsi}'|^2(|\upmu||{\uppsi}'|^2-|\upmu|-\frac{c_0}{r})\leq 0\end{eqnarray}
Manipulating (\ref{3.8}), one obtains
\begin{equation}
\label{3.9}
\frac{1}{2}\mathrm{d}^{\ast}\mathrm{d}|{\uppsi}'|^2-\frac{1}{|\upmu|}<\mathrm{d}|\upmu|, \mathrm{d} |{\uppsi}'|^2>+\frac{r}{2}|\upmu\|{\uppsi}'|^2(|{\uppsi}'|^2-1-\frac{c_1}{r})\leq 0
\end{equation}
where $c_1>0$ is a constant depending on $c_0$. An application of the maximum principle to (\ref{3.9}) yields
\begin{equation}
\label{3.10}
|{\uppsi}'|^2\leq 1+\frac{c_1}{r}
\end{equation}
from which the first bullet of Lemma \ref{l4.1} follows immediately.
\par As for the claimed estimate on the norm of $\upbeta$, start by contracting (\ref{3.6}) first with $(\upalpha,0)$ and then with $(0,\upbeta)$ to get
\begin{eqnarray}
\label{3.11}
\nonumber&&\frac{1}{2}\mathrm{d}^{\ast}\mathrm{d}|\upalpha|^2+|\nabla\upalpha|^2+\frac{r}{2}|\upalpha|^2(|\upalpha|^2+|\upbeta|^2-|\upmu|)+\upkappa_1|\upalpha|^2+\upkappa_2(\upalpha,\upbeta)\\
\nonumber&&+\upkappa_3(\upalpha,\nabla\upalpha)+\upkappa_4(\upalpha,\nabla\upbeta)= 0\\
\nonumber&&\frac{1}{2}\mathrm{d}^{\ast}\mathrm{d}|\upbeta|^2+|\nabla\upbeta|^2+\frac{r}{2}|\upbeta|^2(|\upalpha|^2+|\upbeta|^2+|\upmu|)+\mathrm{\upkappa_1}'(\upbeta,\upalpha)+\mathrm{\upkappa_2}'|\upbeta|^2\\
&&+\mathrm{\upkappa_3}'(\upbeta,\nabla\upalpha)+\mathrm{\upkappa_4}'(\upbeta,\nabla\upbeta)= 0
\end{eqnarray}
where $\mathrm{\upkappa_i}$'s and $\mathrm{\upkappa_i}'$'s depend only on the Riemannian metric. Then, the equations in (\ref{3.11}) yield the following equations in terms of $\upalpha'$ and $\upbeta'$:
\begin{eqnarray}
\label{3.12}
\nonumber&&\frac{1}{2}\mathrm{d}^{\ast}\mathrm{d}|\upalpha'|^2+|\nabla\upalpha'|^2+\frac{r}{2}|\upmu||\upalpha'|^2(|\upalpha'|^2+|\upbeta'|^2-1)+\uplambda_1|\upalpha'|^2\\
\nonumber&&+\uplambda_2(\upalpha',\upbeta')+\uplambda_3(\upalpha',\nabla\upalpha')+\uplambda_4(\upalpha',\nabla\upbeta')= 0\\
\nonumber&&\frac{1}{2}\mathrm{d}^{\ast}\mathrm{d}|\upbeta'|^2+|\nabla\upbeta'|^2+\frac{r}{2}|\upmu||\upbeta'|^2(|\upalpha'|^2+|\upbeta'|^2+1)+\mathrm{\uplambda_1}'(\upbeta',\upalpha')\\
&&+\mathrm{\uplambda_2}'|\upbeta'|^2+\mathrm{\uplambda_3}'(\upbeta',\nabla\upalpha')+\mathrm{\uplambda_4}'(\upbeta',\nabla\upbeta')=0
\end{eqnarray}
where $\mathrm{\uplambda_i}$'s and $\mathrm{\uplambda_i}'$'s depend only on the Riemannian metric.
\par Now, introduce $\mathrm{w}=1-|\upalpha'|^2$. Then, the top equation in (\ref{3.12}) can be rewritten as
\begin{eqnarray}
\label{3.13}
\nonumber-\frac{1}{2}\mathrm{d}^{\ast}\mathrm{d}\mathrm{w}+|\nabla\upalpha'|^2-\frac{r}{2}|\upmu||\upalpha'|^2\mathrm{w}+\frac{r}{2}|\upmu||\upalpha'|^2|\upbeta'|^2+\\
\uplambda_1|\upalpha'|^2+\uplambda_2(\upalpha',\upbeta')+\uplambda_3(\upalpha',\nabla\upalpha')+\uplambda_4(\upalpha',\nabla\upbeta')=0.
\end{eqnarray}
\par Using the estimate in (\ref{3.10}), manipulating the lower order terms and maximizing positive valued functions that do not depend on the value of $r$ or the particular solution $(\upalpha,\upbeta)$, the bottom equation in (\ref{3.12}) and the equation (\ref{3.13}) yield the following inequalities:
\begin{eqnarray}
\label{3.14}
\nonumber&&-\frac{1}{2}\mathrm{d}^{\ast}\mathrm{d}\mathrm{w}+\upzeta_0|\nabla\upalpha'|^2-\frac{r}{2}|\upmu||\upalpha'|^2\mathrm{w}\leq \upzeta_1+\upzeta_2|\nabla\upbeta'|^2\\
\nonumber&&\frac{1}{2}\mathrm{d}^{\ast}\mathrm{d}|\upbeta'|^2+\upeta_0|\nabla\upbeta'|^2+\frac{r}{2}\upeta_1|\upmu||\upbeta'|^2+\frac{r}{2}|\upmu||\upalpha'|^2|\upbeta'|^2\leq \frac{\upeta_2}{r}+\frac{\upeta_3}{r}|\nabla \upalpha'|^2\\
\end{eqnarray}
where $\mathrm{\upzeta_i}$'s and $\mathrm{\upeta_i}$'s are positive constants depending only on the Riemannian metric and the constant $c_0$.
\par Multiplying the top inequality in (\ref{3.14}) by $\frac{\mathrm{k}}{r}$ where $\mathrm{k}$ is a positive constant large enough to satisfy
\begin{itemize}
\item $\mathrm{k}\upzeta_0\geq\upeta_3$ and
\item $\upeta_0\geq\mathrm{k}\upzeta_2$,
\end{itemize}
and adding the resulting inequality to the bottom inequality in (\ref{3.14}), we deduce that there are positive constants $c_2$ and $c_3$ that depend only on the Riemannian metric and the constant $c_0$ such that
\begin{equation}
\label{3.15}
\mathrm{d}^{\ast}\mathrm{d}(|\upbeta'|^2-\frac{c_2}{r}\mathrm{w}-\frac{c_3}{r^2})+r|\upmu||\upalpha'|^2(|\upbeta'|^2-\frac{c_2}{r}\mathrm{w}-\frac{c_3}{r^2})\leq0.
\end{equation}
Then, an application of the maximum principle to (\ref{3.15}) yields
$$|\upbeta'|^2\leq\frac{{c_2}}{r}(1-|\upalpha'|^2)+\frac{{c_3}}{r^2}$$
which, eventually, gives rise to the second bullet of Lemma \ref{l4.1} after multiplying both sides of the inequality by $|\upmu|$. \QED
\par Given Lemma \ref{l4.1}, the next lemma finds a priori bounds on the derivatives of $\upalpha$ and $\upbeta$.
\begin{lemma} \label{l4.2}Fix a bound on the $\mathrm{C}^3$-norm of $\upmu$. Given $r\geq 1$ and $t\in \mathrm{S^1}$, let $(\mathbb{A},\uppsi=(\upalpha,\upbeta))$ denote a solution of the $t$ and $r$ version of the equations in (\ref{3.1}). Then, for each integer $n\geq1$ there exists a constant $c_n\geq 1$, which is independent of the value of $t\in \mathrm{S^1}$, the value of $r\geq 1$ and the solution $(\mathbb{A},\uppsi=(\upalpha,\upbeta))$, with the following significance:
\begin{itemize}
\item $|\nabla^n \upalpha|\leq c_n r^{n/2}$
\item $|\nabla^n \upbeta|\leq c_n r^{(n-1)/2}.$
\end{itemize}
The following is also true: Fix $\upepsilon>0$. There exists $\updelta>0$ and $\upkappa>1$ such that if $r>\upkappa$ and if $|\upalpha|\geq |\upmu|^{1/2}-\updelta$ in any given ball of radius $2\upkappa r^{-1/2}$ in $\mathrm{M}_t$, then $|\nabla^n \upalpha|\leq \upepsilon c_n r^{n/2}$ for $n\geq 1$ and $|\nabla^n \upbeta|\leq \upepsilon c_n r^{(n-1)/2}$ for all $n\geq0$ in the concentric ball with radius $\upkappa r^{-1/2}$.
\end{lemma}
\par\noindent\textbf{Proof.} The proof is essentially identical to that of Lemma~2.3 in \cite{t1}. This is to say that 
the proof is local in nature:  Fix a Gaussian coordinate chart centered at any given point 
in $\mathrm{M}$ so as to view the equations in (\ref{3.1}) as equations on a small ball in $\mathbb{R}^3$.  Then rescale 
coordinates by writing $x = r^{-1/2}y$ so that the resulting equations are on a ball of radius 
$\mathcal{O}(r^{1/2})$ in $\mathbb{R}^3$.  The $r$-dependence of these rescaled equations is such that standard elliptic 
regularity techniques provide uniform bounds on the rescaled versions of $\upbeta$ and the 
derivatives of the rescaled $\upalpha$ and $\upbeta$ in the unit radius ball about the origin.  Rescaling back to the original coordinates will give what is 
claimed by the lemma.  
\QED

\par One of the key implications of Lemma \ref{l4.1} is a priori bounds on the values of $\mathcal{E}$. First, note that since $\upnu\wedge\upmu>0$ at each $t\in \mathrm{S^1}$, it follows that 
\begin{equation}\label{4.41}
\upnu=\ast\frac{\mathrm{q}}{|\upmu|}\upmu+\upupsilon 
\end{equation}
where $\mathrm{q}=<\upnu,\ast\upmu>|\upmu|^{-1}$ is a positive valued function on $\mathrm{M}_t$ at each $t\in \mathrm{S^1}$, and $\upupsilon\wedge\upmu=0$. We use (\ref{4.41}) in the following proof of Lemmas \ref{p3.5}.
\\
\par\noindent\textbf{Proof of Lemma \ref{p3.5}.} Fix $r\geq1$ and $t\in\mathrm{S^1}$. Let $(\mathbb{A},\uppsi)$ be a solution of the $t$ and $r$ version of the equations in (\ref{3.1}). Write $\mathbb{A}=\mathbb{A_S}+2\mathrm{a}$ and $\uppsi=(\upalpha,\upbeta)$. Then, by (\ref{4.41}) we can write
\begin{equation}\label{4.42}
\mathcal{E}(\mathbb{A},\uppsi)=i\int_{\mathrm{M}}\upnu\wedge\mathrm{da}=r\int_{\mathrm{M}}\mathrm{q}(|\upmu|-|\upalpha|^2)+i\int_{\mathrm{M}}\upupsilon\wedge\mathrm{da}.
\end{equation}
Now, it follows from (\ref{3.1}) and Lemma \ref{l4.1} that \begin{equation}
\label{3.21}
\mathcal{E}(\mathbb{A},\uppsi)\geq \frac{1}{2}r\int_{\mathrm{M}}\mathrm{q}(|\upmu|-|\upalpha|^2)-c_4\geq-c_5
\end{equation}
where $c_4,c_5>0$ are constants depending only on the Riemannian metric. 
\QED
\par \noindent\textbf{Proof of Propostions \ref{p3.1} and \ref{p3.7}.} Proposition \ref{p3.1} follows directly from Proposition \ref{p3.7}. Given Lemma \ref{l4.1}, the proof of the latter is identical but for minor changes to the proof of Theorem 2.1 given in Section 6 of \cite{t1}. The proof of the second bullet is proved just as in Lemma 6.5 in \cite{t1}.\QED

\par \noindent\textbf{Proof of Proposition \ref{p3.2}.} In the case when $c_1(\mathrm{E})\neq0$, the claim about $|\uppsi|$ follows from Lemma \ref{l4.1} given that $\upalpha$ is a section of $\mathrm{E}$. This understood, we now assume that $\mathrm{E}=\underline{\mathbb{C}}$. To start, let $1_{\underline{\mathbb{C}}}$ denote a unit length trivializing section of the $\underline{\mathbb{C}}$ summand. There exists a unique connection $\mathrm{A}_0$ on $\mathrm{K^{-1}}$ such that the section $\uppsi_0=(1_{\underline{\mathbb{C}}},0)$ of $\mathbb{S}_0=\underline{\mathbb{C}}\oplus\mathrm{K}^{-1}$ obeys $\mathcal{D}_{\mathrm{A}_0}\uppsi_0=0$. Now, we look for a solution of the equations in (\ref{3.1}) of the form 

$$(\mathbb{A},\uppsi)=(\mathrm{A_0}+2(2r)^{1/2}\mathrm{b},|\upmu|^{1/2}\uppsi_0+\upphi)$$ 

\noindent with $(\mathrm{b},\upphi)\in\mathrm{C^{\infty}(M};i\mathrm{T^*M}\oplus\mathbb{S})$. Then, $(\mathbb{A}, \uppsi)$ will solve the equations in (\ref{3.1}) if $\mathfrak{b}=(\mathrm{b},\upphi,g)\in\mathrm{C^{\infty}(M};i\mathrm{T^*M}\oplus\mathbb{S}\oplus i\underline{\mathbb{R}})$ solves the following system of equations:
\begin{align}\label{5.1}
\nonumber\ast\mathrm{db}-\mathrm{d}g-2^{-1/2}r^{1/2}[|\upmu|^{1/2}({\uppsi_0}^{\dagger}\uptau\upphi+\upphi^{\dagger}\uptau\uppsi_0)+\upphi^{\dagger}\uptau\upphi]&=-2^{-3/2}r^{-1/2}\ast\mathrm{F_{A_0}}\\
\nonumber\mathcal{D}_{\mathrm{A}_0}\upphi+2^{1/2}r^{1/2}[|\upmu|^{1/2}(\mathfrak{cl}(\mathrm{b})\uppsi_0+g\uppsi_0)+(\mathfrak{cl}(\mathrm{b})\upphi+g\upphi)]&=-\mathfrak{cl}(\mathrm{d|\upmu|^{1/2}})\uppsi_0\\
\nonumber-\mathrm{d^{\ast}b}-2^{-1/2}|\upmu|^{1/2}r^{1/2}(\upphi^{\dagger}\uppsi_0-{\uppsi_0}^{\dagger}\upphi)&=0.\\
\end{align}
For notational convenience, we denote by $\mathcal{L}_0$ the operator $\mathcal{L}_{(\mathrm{A}_0,|\upmu|^{1/2}\uppsi_0)}$ as defined in (\ref{4.18}). Then, the equations in (\ref{5.1}) can be rewritten as
\begin{eqnarray}\label{5.2}
\nonumber\mathcal{L}_0(\mathrm{b},\upphi,g)+r^{1/2}\left(\begin{array}{c}-2^{-1/2}\upphi^{\dagger}\uptau\upphi\\ 2^{1/2}(\mathfrak{cl}(\mathrm{b})\upphi+g\upphi)\\
0\end{array}\right)=\left(\begin{array}{c}-2^{-3/2}r^{-1/2}\ast\mathrm{F_{A_0}}\\ -\mathfrak{cl}(\mathrm{d|\upmu|^{1/2}})\uppsi_0\\
0\end{array}\right).\\
\end{eqnarray}
Now, for $\mathfrak{b}=(\mathrm{b},\upphi,g)$ and $\mathfrak{b}'=(\mathrm{b}',\upphi',g')$ in $\mathrm{C^{\infty}(M};i\mathrm{T^*M}\oplus\mathbb{S}\oplus i\underline{\mathbb{R}})$, let $(\mathfrak{b},\mathfrak{b}')\mapsto\mathfrak{b}\ast\mathfrak{b}'$ be the bilinear map defined by
\begin{equation}\label{5.3}
\mathfrak{b}\ast\mathfrak{b}'=\frac{1}{2}\left(\begin{array}{c}-2^{-1/2}(\upphi^{\dagger}\uptau\upphi'+\upphi'^{\dagger}\uptau\upphi)\\
2^{1/2}(\mathfrak{cl}(\mathrm{b})\upphi'+g\upphi'+\mathfrak{cl}(\mathrm{b'})\upphi+g'\upphi)\\
0\end{array}\right),
\end{equation}
and let $\mathfrak{u}$ denote the section defined by $(-2^{-3/2}r^{-1/2}\ast\mathrm{F_{A_0}},-\mathfrak{cl}(\mathrm{d|\upmu|^{1/2}})\uppsi_0,0)$ of $i\mathrm{T^*M}\oplus\mathbb{S}\oplus i\underline{\mathbb{R}}$. Then, (\ref{5.2}) has the schematic form 
\begin{equation}\label{schematic}
\mathcal{L}_0\mathfrak{b}+r^{1/2}\mathfrak{b}\ast\mathfrak{b}=\mathfrak{u}.
\end{equation}  
Our plan is to use the contraction mapping theorem to solve (\ref{schematic}) in a manner much like what is done in the proof of Proposition 2.8 of \cite{t2}. To set the stage for this, we first introduce the Hilbert space $\mathbb{H}$ as the completion of $\mathrm{C^{\infty}(M};i\mathrm{T^*M}\oplus\mathbb{S}\oplus i\underline{\mathbb{R}})$ with respect to the norm whose square is:
\begin{equation}\label{5.4}
{||\upxi||_{\mathbb{H}}}^2=\int_{\mathrm{M}}|\nabla_0\upxi|^2+\frac{1}{4}r\int_{\mathrm{M}}|\upxi|^2,
\end{equation}
where $\nabla_0$ denotes the covariant derivative on sections of $i\mathrm{T^*M}\oplus\mathbb{S}\oplus i\underline{\mathbb{R}}$ that acts as the Levi-Civita covariant derivative on sections of $i\mathrm{T^*M}$, the covariant derivative defined by $\mathrm{A_0}$ on sections of $\mathbb{S}$, and that defined by the exterior derivative on sections of $i\underline{\mathbb{R}}$.

\begin{lemma}\label{l4.3}
There exists $\upkappa\geq 1$ such that
\begin{itemize}
\item $||\upxi||_6\leq \upkappa ||\upxi||_{\mathbb{H}}$ and $||\upxi||_4\leq\upkappa r^{-1/8}||\upxi||_{\mathbb{H}}$ for all $\upxi\in\mathbb{H}$.
\item If $r\geq\upkappa$, then $\upkappa^{-1}||\upxi||_{\mathbb{H}}\leq||\mathcal{L}_0\upxi||_2\leq\upkappa||\upxi||_{\mathbb{H}}$ for all $\upxi\in\mathbb{H}$.
\end{itemize}
\end{lemma}
\par\noindent\textbf{Proof.} The first bullet follows using a standard Sobolev inequality with the fact that $|\mathrm{d}|\upxi||\leq|\nabla_0\upxi|$. The right hand inequality in the second bullet follows by simply from the appearance of only first derivatives in $\mathcal{L}_0$. To obtain the left hand inequality of the second bullet, use the Bochner-type formula for the operator ${\mathcal{L}_0}^2$ (see (5.21) in \cite{t2}). To elaborate, let $\mathrm{f}$ be any given function on $\mathrm{M}$. Write a section $\upxi$ of $i\mathrm{T^*M}\oplus\mathbb{S}\oplus i\underline{\mathbb{R}}$ as $(\mathrm{b},\upphi,g)$. Then, ${\mathcal{L}_{(\mathrm{A}_0,\mathrm{f}\uppsi_0)}}^2(\mathrm{b},\upphi,g)$ has respective $i\mathrm{T^*M}$, $\mathbb{S}$ and $i\underline{\mathbb{R}}$ components
\begin{eqnarray}\label{4.20}
\nonumber&&\nabla^{\dagger}\nabla\mathrm{b}+\textcolor{blue}{2r\mathrm{f^2}\mathrm{b}}+r^{1/2}\mathbb{V}_1(\upxi)\\
\nonumber&&{\nabla_{\mathrm{A_0}}}^{\dagger}{\nabla_{\mathrm{A_0}}}\upphi+\textcolor{blue}{2r\mathrm{f^2}\upphi}+r^{1/2}\mathbb{V}_2(\upxi)\\
&&\mathrm{d^{\ast}d}g+\textcolor{blue}{2r\mathrm{f^2}g}+r^{1/2}\mathbb{V}_3(\upxi),
\end{eqnarray}
where $\mathbb{V}_{\mathrm{i}}$ are zero'th order endomorphisms with absolute value bounded by an $r$-independent constant. In the case at hand, $\mathrm{f}=|\upmu|^{1/2}$ is strictly bounded away from zero. This last point understood, then the left hand inequality in the second bullet of the lemma follows by first taking the $\mathrm{L}^2$ inner product of ${\mathcal{L}_0}^2\upxi$ with $\upxi$ and then integrating  by parts to rewrite the resulting integral.
\QED
\par It follows from Lemma \ref{l4.3} that the operator $\mathcal{L}_0$ is invertible when $r$ is large. This understood, write $\mathfrak{y}={\mathcal{L}_0}^{-1}\mathfrak{u}$,

\begin{lemma}\label{l4.4}
There exists $\upkappa\geq 1$ for use in Lemma \ref{l4.3} such that when $r\geq\upkappa$, then the corresponding $\mathfrak{y}={\mathcal{L}_0}^{-1}\mathfrak{u}$ obeys $|\mathfrak{y}|\leq c_0r^{-1/2}$.
\end{lemma}
\noindent\textbf{Proof.} Let $\Delta$ denote the operator that is obtained from what is written in the $\mathrm{f}=|\upmu|^{1/2}$ version of (\ref{4.20}) by setting $\mathbb{V}_{\mathrm{i}}$ all equal to zero. The latter has Green's function $\mathrm{G}$, a positive, symmetric function on $\mathrm{M\times M}$ with pole along the diagonal. Moreover, there exists an $r$-independent constant $c>1$ such that if $\mathrm{x,y}\in\mathrm{M}$, then
\begin{eqnarray}\label{4.50}
\nonumber &&\mathrm{G}(\mathrm{x},\mathrm{y})\leq \frac{c}{dist(\mathrm{x},\mathrm{y})}e^{-\sqrt{r}\frac{dist(\mathrm{x},\mathrm{y})}{c}},\\
&&|\mathrm{dG}|(\mathrm{x},\mathrm{y})\leq c(\frac{1}{dist(\mathrm{x},\mathrm{y})^{2}}+\frac{\sqrt{r}}{dist(\mathrm{x},\mathrm{y})})e^{-\sqrt{r}\frac{dist(\mathrm{x},\mathrm{y})}{c}}.
\end{eqnarray}
Both of these bounds follow by using the maximum principle with a standard parametrix for $\mathrm{G}$ near the diagonal in $\mathrm{M\times M}$.
\par Now write (\ref{4.20}) as $\Delta\upxi+r^{1/2}\mathbb{V}\upxi$, and then use $\mathrm{G}$, the fact that ${\mathcal{L}_0}^2\mathfrak{y}=\mathcal{L}_0\mathfrak{u}$, and the uniform bounds on the terms $\mathbb{V}_{\mathrm{i}}$ to see that
\begin{equation}
\nonumber |\mathfrak{y}|(\mathrm{x})\leq c'\int_{\mathrm{M}}\mathrm{G}(\mathrm{x},\cdot)(1+r^{1/2}(1+|\mathfrak{y}|)),
\end{equation}
where $c'$ is independent of $r$. This last equation together with (\ref{4.50}) yields
\begin{equation}
\nonumber |\mathfrak{y}|(\mathrm{x})\leq c''r^{-1/2}(1+\mathrm{sup_M}|\mathfrak{y}|),
\end{equation}
where $c''$ is also independent of $r$. The lemma follows from this bound.\QED

\par With $\mathfrak{y}$ in hand, it follows that $\upxi\in\mathbb{H}$ is a solution of the equations in (\ref{schematic}) if $\tilde{\upxi}=\upxi-\mathfrak{y}$ is a solution of the equation $\mathcal{L}_0\mathfrak{\tilde{\upxi}}+r^{1/2}(\tilde{\upxi}\ast\tilde{\upxi}+2\mathfrak{y}\ast\tilde{\upxi})=-r^{1/2}\mathfrak{y}\ast\mathfrak{y}$. To find a solution $\tilde{\upxi}$ of the latter equation, introduce the map $\mathbb{T:H\rightarrow H}$ defined by
\begin{equation}\label{5.7}
\mathbb{T}:\mathfrak{\tilde{\upxi}}\mapsto-r^{1/2}{\mathcal{L}_0}^{-1}(\mathfrak{y}\ast\mathfrak{y}+\tilde{\upxi}\ast\tilde{\upxi}+2\mathfrak{y}\ast\tilde{\upxi}).
\end{equation}
Note in this regard that Sobolev inequalities in Lemma \ref{l4.3} guarantee that $\mathbb{T}$ does indeed define a smooth map from $\mathbb{H}$ onto itself when $r$ is larger than some fixed constant. Our goal now is to show that the map $\mathbb{T}$ has a unique fixed point with small norm.
Given $\mathrm{R}\geq1$, we let $\mathrm{B_R}\in\mathbb{H}$ denote the ball of radius $r^{-1/2}\mathrm{R}$ centered at the origin. We next invoke 
\begin{lemma}\label{l4.5}
There exists $\upkappa>1$, and given $\mathrm{R}\geq\upkappa$, there exists $\upkappa_{\mathrm{R}}$ such that if $r\geq\upkappa_{\mathrm{R}}$, then $\mathbb{T}$ maps $\mathrm{B_{\mathrm{R}}}$ onto itself as a contraction mapping.
\end{lemma}
\par\noindent\textbf{Proof.} Let $\mathrm{R}>1$ be such that $||\mathfrak{y}||_{\infty}\leq\frac{1}{2^{10}} r^{-1/2}\mathrm{R}^{1/2}$. We first show that if $r$ is large, then $\mathbb{T}$ maps $\mathrm{B_{\mathrm{R}}}$ into itself. Indeed, this follows from Lemma \ref{l4.3} using the following chain of inequalities:
\begin{align}\label{5.8}
\nonumber||\mathbb{T}(\tilde{\upxi})||_{\mathbb{H}}&\leq ||-r^{1/2}\mathfrak{y}\ast\mathfrak{y}-r^{1/2}(\tilde{\upxi}\ast\tilde{\upxi}+2\mathfrak{y}\ast\tilde{\upxi})||_2\\
\nonumber&\leq r^{1/2}||\mathfrak{y}\ast\mathfrak{y}||_2+r^{1/2}||\tilde{\upxi}\ast\tilde{\upxi}+2\mathfrak{y}\ast\tilde{\upxi}||_2\\
\nonumber&\leq \frac{1}{4}r^{-1/2}\mathrm{R}+r^{1/2}(||\tilde{\upxi}\ast\tilde{\upxi}||_2+2||\mathfrak{y}\ast\tilde{\upxi}||_2)\\
\nonumber&\leq \frac{1}{4}r^{-1/2}\mathrm{R}+r^{1/2}({||\tilde{\upxi}||_4}^2+2||\mathfrak{y}||_4||\tilde{\upxi}||_4)\\
\nonumber&\leq \frac{1}{4}r^{-1/2}\mathrm{R}+r^{1/2}(\upkappa r^{-1/4}{||\tilde{\upxi}||_{\mathbb{H}}}^2+r^{-1/2}\mathrm{R}^{1/2}\upkappa r^{-1/8}||\tilde{\upxi}||_{\mathbb{H}})\\
\nonumber&\leq \frac{1}{4}r^{-1/2}\mathrm{R}+r^{1/2}(\upkappa r^{-1/4}r^{-1}\mathrm{R}^2+r^{-1/2}\mathrm{R}^{1/2}\upkappa r^{-1/8}r^{-1/2}\mathrm{R})\\
&\leq r^{-1/2}\mathrm{R}(\frac{1}{4}+2\upkappa \mathrm{R}r^{-1/8}).
\end{align}
Next, using similar arguments, we show that $\mathbb{T}|_{\mathrm{B_R}}$ is a contraction mapping. In this regard, let $\tilde{\upxi}_1, \tilde{\upxi}_2\in\mathrm{B_R}$, then
\begin{align}\label{5.9}
\nonumber||\mathbb{T}(\tilde{\upxi}_1)-\mathbb{T}(\tilde{\upxi}_2)||_{\mathbb{H}}&\leq||-r^{1/2}(\tilde{\upxi}_1\ast\tilde{\upxi}_1+2\mathfrak{y}\ast\tilde{\upxi}_1)+r^{1/2}(\tilde{\upxi}_2\ast\tilde{\upxi}_2+2\mathfrak{y}\ast\tilde{\upxi}_2)||_2\\
\nonumber&\leq r^{1/2}(||(\tilde{\upxi}_1\ast\tilde{\upxi}_1-\tilde{\upxi}_2\ast\tilde{\upxi}_2)||_2+2||\mathfrak{y}\ast\tilde{\upxi}_1-\mathfrak{y}\ast\tilde{\upxi}_2||_2)\\
\nonumber&\leq r^{1/2}(||(\tilde{\upxi}_1+\tilde{\upxi}_2)\ast(\tilde{\upxi}_1-\tilde{\upxi}_2)||_2+||\mathfrak{y}\ast(\tilde{\upxi}_1-\tilde{\upxi}_2)||_2)\\
\nonumber&\leq r^{1/2}(||\tilde{\upxi}_1+\tilde{\upxi}_2||_4||\tilde{\upxi}_1-\tilde{\upxi}_2||_4+2||\mathfrak{y}||_4||\tilde{\upxi}_1-\tilde{\upxi}_2||_4)\\
\nonumber&\leq r^{1/2}(||\tilde{\upxi}_1||_4+||\tilde{\upxi}_2||_4+2||\mathfrak{y}||_4)||\tilde{\upxi}_1-\tilde{\upxi}_2||_4\\
\nonumber&\leq r^{1/2}(2\upkappa r^{-1/8}r^{-1/2}\mathrm{R}+r^{-1/2}\mathrm{R}^{1/2})\upkappa r^{-1/8}||\tilde{\upxi}_1-\tilde{\upxi}_2||_{\mathbb{H}}\\
&\leq 3\upkappa^2\mathrm{R}r^{-1/8}||\tilde{\upxi}_1-\tilde{\upxi}_2||_{\mathbb{H}}.
\end{align}
Therefore, by the contraction mapping theorem, there exists a unique fixed point of the map $\mathbb{T}$ in the ball $\mathrm{B_R}$. Moreover, by standard elliptic regularity arguments, it follows that the fixed point is smooth, therefore it is an element of $\mathrm{C^{\infty}(M};i\mathrm{T^*M}\oplus\mathbb{S}\oplus i\underline{\mathbb{R}})$.\QED
\par We next find an $r$-independent constant $\upkappa$ and prove that the norm of $\uppsi=|\upmu|^{1/2}\uppsi_0+\upphi$ is bounded from below by $|\upmu|^{1/2}-\upkappa r^{-1/2}$. To this end, note that $\tilde{\upxi}$ obeys the equation
\begin{equation}\label{4.51}
\Delta\tilde{\upxi}+r^{1/2}\mathbb{V}\tilde{\upxi}=-r^{1/2}\mathcal{L}_0(\mathfrak{y}\ast\mathfrak{y}+\tilde{\upxi}\ast\tilde{\upxi}+2\mathfrak{y}\ast\tilde{\upxi}).
\end{equation}
What with (\ref{4.50}) and the bound $|\mathfrak{y}|\leq2r^{-1/2}\mathrm{R}$ this last equation implies is 
\begin{align}\label{4.52}
\nonumber|\tilde{\upxi}|(\mathrm{x})\leq c_0r^{-1/2}+c_0r^{1/2}\int_{\mathrm{M}}(\frac{1}{dist(\mathrm{x},\cdot)^2}+\frac{\sqrt{r}}{dist(\mathrm{x},\cdot)})e^{-\sqrt{r}\frac{dist(\mathrm{x},\cdot)}{c}}(|\tilde{\upxi}|^2+r^{-1/2}|\tilde{\upxi}|)]\\
\end{align}
where $c_0$ is independent of $x$ and $r$. Bound the term $r^{-1/2}|\tilde{\upxi}|$ in the integral by $|\tilde{\upxi}|^2+r^{-1}$. The contribution to the right hand side of (\ref{4.52}) of the resulting term with $r^{-1}$ factor is bounded by $c_1r^{-1/2}$ where $c_1$ is independent of $r$. To say something about the term with $|\tilde{\upxi}|^2$, note that the function $\frac{1}{dist(\mathrm{x},\cdot)}|\tilde{\upxi}|$ is square integrable with $\mathrm{L^2}$-norm bounded by an $x$-independent multiple of the ${\mathrm{L}^2}_1$-norm of $|\tilde{\upxi}|$; and thus by $c_2||\tilde{\upxi}||_{\mathbb{H}}$ with $c_2$ independent of $r$ and $\tilde{\upxi}$. This understood, the term in the integral with $|\tilde{\upxi}|^2$ contributes at most $c_3(r^{1/2}{||\tilde{\upxi}||_{\mathbb{H}}}^2+r||\tilde{\upxi}||_2||\tilde{\upxi}||_{\mathbb{H}})$ with $c_3$ independent of $r$ and $\tilde{\upxi}$. The latter is bounded by an $r$-independent multiple of $r^{-1/2}$. Thus, we see that $|\tilde{\upxi}|\leq c_4 r^{-1/2}$ which proves our claim that $|\uppsi|\geq|\upmu|^{1/2}-\upkappa r^{-1/2}$.
\par We now turn to the claim about uniqueness. To this end, let $\updelta\in(0,\frac{\mathrm{inf}_{\mathrm{M}}|\upmu|}{2})$ and let $(\mathbb{A},\uppsi)$ be a solution of some $t\in\mathrm{S}^1$ and $r\geq1$ version of the equations in (\ref{3.1}) with the property that $|\uppsi|\geq|\upmu|^{1/2}-\updelta$ at each point in $\mathrm{M}$. Granted such is the case, it follows from Lemma \ref{l4.1} that $|\upalpha|\geq|\upmu|^{1/2}-\updelta-\upkappa r^{-1/2}$ at each point in $\mathrm{M}$, with $C_0$ independent of $r$. We now make use of Lemma \ref{l4.2} to see the following: Given $\upepsilon>0$, there exists $\updelta_{\upepsilon}>0$ such that if $\updelta<\updelta_{\upepsilon}$, then 
\begin{eqnarray}\label{5.16}
\nonumber&&|\upmu|^{1/2}-\upepsilon\leq|\upalpha|\leq|\upmu|^{1/2}+\upepsilon\;\mathrm{and}\;|\upbeta|\leq\upepsilon r^{-1/2},\\
\nonumber&&|\nabla\upalpha|\leq\upepsilon r^{1/2}\;\mathrm{and}\;|\nabla\upbeta|\leq\upepsilon,\\
&&|\nabla^2\upalpha|\leq\upepsilon r\;\mathrm{and}\;|\nabla^2\upbeta|\leq\upepsilon r^{1/2}.
\end{eqnarray}
\par Since $\upalpha$ is nowhere zero for sufficiently large $r>1$, one has $u=\bar{\upalpha}/|\upalpha|\in \mathcal{G}$. Now, change $(\mathbb{A},\uppsi)$ to a new gauge by $u$, and denote the resulting pair of gauge and spinor fields again by $(\mathbb{A},\uppsi)$. Since $u\upalpha=|\upalpha|1_{\underline{\mathbb{C}}}$, one has $\mathbb{A}=\mathrm{A}_0+2i\mathrm{a}$ where
\begin{equation}\label{5.17}
\mathrm{a}=-\frac{i}{2}(\upalpha^{-1}\nabla\upalpha-\bar{\upalpha}^{-1}\nabla\bar{\upalpha}).
\end{equation}
Then, (\ref{5.16}) and (\ref{5.17}) imply
\begin{equation}\label{5.18}
r^{-1/2}|\mathrm{a}|+r^{-1}|\nabla\mathrm{a}|\leq c_0\upepsilon.
\end{equation}
\par We now change $(\mathbb{A},\uppsi)$ to yet another gauge so as to write the resulting pair of connection and spinor as $(\mathrm{A}_0+2(2r)^{1/2}\mathrm{b},|\upmu|^{1/2}\uppsi_0+\upphi)$ where $(\mathrm{b},\upphi,0)$ obey (\ref{5.1}). This gauge transformation is written $e^{i\mathrm{x}}$ where $\mathrm{x}:\mathrm{M}\rightarrow\mathbb{R}$. Thus, the pair $(\mathrm{b},\upphi)$ is 
\begin{eqnarray}\label{5.19}
\nonumber\mathrm{b}&=&i(2r)^{-1/2}(\mathrm{a}-\mathrm{dx})\\
\upphi&=&e^{i\mathrm{x}}\uppsi-|\upmu|^{1/2}\uppsi_0.
\end{eqnarray}
Equation  (\ref{5.1})  is obeyed if and only if $\mathrm{x}$ obeys the equation
\begin{equation}\label{5.20}
\mathrm{d}^{\ast}\mathrm{dx}+2|\upmu|^{1/2}r|\upalpha|\sin{\mathrm{x}}=\mathrm{d^{\ast}b}.
\end{equation}
We can now proceed along the lines of what is done in \cite{t2} to solve an analogous equation, namely (2.16) in \cite{t2}. In particular, the arguments in \cite{t2} can be used with only small modifications to find an $r$-independent constant $\upkappa$ such that if the constant $\upepsilon$ in (\ref{5.16}) is bounded by $\upkappa^{-1}$ and $r\geq \upkappa$, then (\ref{5.20}) has a unique solution, $\mathrm{x}$, with 
\begin{equation}\label{5.25}|\mathrm{x}|+r^{1/2}|\mathrm{dx}|\leq \upkappa\upepsilon.\end{equation}
Granted this, it follows that $\mathfrak{b}=(\mathrm{b},\upphi,0)$ with $(\mathrm{b},\upphi)$ as in (\ref{5.19}) obeys (\ref{schematic}) and that 
\begin{equation}\label{4.54}
|\mathfrak{b}|\leq c\upepsilon
\end{equation}
with $c>0$ a constant that is independent of $\upepsilon$ and $r$. Then, $\mathfrak{h}=\mathfrak{b}-\mathfrak{y}$ obeys $\mathcal{L}_0\mathfrak{\mathfrak{h}}=r^{1/2}(\mathfrak{y}\ast\mathfrak{y}+\mathfrak{h}\ast\mathfrak{h}+2\mathfrak{y}\ast\mathfrak{h})$ and $||\mathfrak{h}||_{\infty}\leq c_0\upepsilon$ where $c_0$ is independent of $(\mathbb{A},\uppsi)$ and $r$. This understood, it follows from Lemma \ref{l4.3} that
\begin{equation}\label{4.56}
||\mathfrak{h}||_{\mathbb{H}}\leq\frac{1}{4}\mathrm{R}_{\mathfrak{y}}r^{-1/2}+c_1r^{1/2}||\mathfrak{h}||_{\infty}||\mathfrak{h}||_2\leq \frac{1}{4}\mathrm{R}_{\mathfrak{y}}r^{-1/2}+c_2r^{1/2}\upepsilon||\mathfrak{h}||_2,
\end{equation} 
where $\mathrm{R}_{\mathfrak{y}}$ is an $r$ independent constant such that $||\mathfrak{y}||_{\infty}\leq\frac{1}{2^{10}}r^{-1/2}\mathrm{R}_{\mathfrak{y}}$ and $c_1,c_2>0$ are constants which are both independent of $(\mathbb{A},\uppsi)$ and $r$. This last inequality implies that $||\mathfrak{h}||_{\mathbb{H}}<\mathrm{R}_{\mathfrak{y}}r^{-1/2}$ when $\upepsilon<c_4$ with $c_4$ an $r$ and $(\mathbb{A},\uppsi)$ independent constant. This understood, it follows from Lemma \ref{l4.5} that $(\mathbb{A},\uppsi)$ is gauge equivalent to the solution of (\ref{3.1}) that was constructed from Lemma \ref{l4.5}'s fixed point of the map $\mathbb{T}$ when $r$ is larger than some fixed constant. This then proves the uniqueness assertion made by Proposition \ref{p3.2}.
\par We introduce $(\mathbb{A}_{\underline{\mathbb{C}}},\uppsi_{\underline{\mathbb{C}}})$ to denote the solution that is obtained from Lemma \ref{l4.5}'s fixed point. This solution is of the form $(\mathrm{A_0}+2(2r)^{1/2}\mathrm{b},|\upmu|^{1/2}\uppsi_0+\upphi)$. Our final task is to prove that the $(\mathbb{A}_{\underline{\mathbb{C}}},\uppsi_{\underline{\mathbb{C}}})$ version of the operator in (\ref{4.18}) has trivial kernel. To see that such is the case, remember that $(\mathrm{b},\upphi)$ has norm bounded by $c_0r^{-1/2}$ with $c_0$ independent of $r$. This being the case, the operator in question differs from the operator $\mathcal{L}_0$ by a zero'th order term with bound independent of $r$. As a consequence, there is a constant $c>0$ which is independent of $r$ and such that
\begin{equation}\label{4.60}
||\mathcal{L}_{(\mathbb{A}_{\underline{\mathbb{C}}},\uppsi_{\underline{\mathbb{C}}})}\upxi||_2\geq c||\upxi||_{\mathbb{H}}
\end{equation}
for all $\upxi\in\mathbb{H}$ when $r$ is large. This understood, the fact that $(\mathbb{A}_{\underline{\mathbb{C}}},\uppsi_{\underline{\mathbb{C}}})$ is non-degenerate when $r$ is large follows from Lemma \ref{l4.3}.\QED
\section{Proof of the Main Theorem}
\label{S5}
We prove Proposition \ref{p3.9} in this section and thus complete the proof of our main theorem. The proof that follows has nine parts.\\
\par \emph{Part 1}: Here we say more about the solution of each $t\in\mathrm{S^1}$ version of the equations in (\ref{3.1}) provided by Proposition \ref{p3.2}. We denote this solution as $(\mathbb{A}_{\underline{\mathbb{C}}},\uppsi_{\underline{\mathbb{C}}})$ and write it at times as $(\mathbb{A}_{\underline{\mathbb{C}}}=\mathbb{A}_{\mathbb{S}_0}+2\Ac,\uppsi_{\underline{\mathbb{C}}}=(\upalpha_{\underline{\mathbb{C}}},\upbeta_{\underline{\mathbb{C}}}))$ where $\mathbb{A}_{\mathbb{S}_0}$ is a $t$-independent connection on the line bundle $\mathrm{K}^{-1}=det(\mathbb{S}_0)$ with harmonic curvature form, and where $\mathrm{A}_{\underline{\mathbb{C}}}$ is a connection on the trivial bundle $\underline{\mathbb{C}}$. Since each $t\in\mathrm{S^1}$ version of these solutions is non-degenerate, the family parametrized by $t\in\mathrm{S^1}$ can be changed by $t$-dependent gauge transformations to define a smooth map from the universal cover, $\mathbb{R}$, of $\mathrm{S^1}$ into $\mathcal{C}$. Moreover, because $\upalpha_{\underline{\mathbb{C}}}$ is nowhere zero, a further gauge transformation can be applied if necessary to obtain a $2\pi$-periodic map from $\mathbb{R}$ into $\mathcal{C}$ and thus a map from $\mathrm{S^1}$ into $\mathcal{C}$. This understood, we can view $\Ac$ as a connection on the trivial bundle over $\mathrm{S^1\times M}$. We write its curvature form as 
\begin{equation}\label{5.51}
\FAc=\FAct+\mathrm{dt}\wedge\dAc.
\end{equation}
where $\FAct$ denotes the component long $\mathrm{M}_t$.
Note that the integral of $\frac{i}{2\pi}\upomega\wedge\mathrm{dt}\wedge\dAc$ over $\mathrm{S^1\times M}$ is zero since $(\mathbb{A}_{\underline{\mathbb{C}}},\uppsi_{\underline{\mathbb{C}}})$ is a 1-parameter family of solutions of the equations in (\ref{3.1}). To see this, use an integration by parts, the fact that $\mathrm{d}\upnu=\dot{\upmu}$ and the equation in (\ref{4.3}) to get

\begin{eqnarray}\label{5.52}
\nonumber \frac{i}{2\pi}\int_{\mathrm{S^1\times M}}\upomega\wedge\mathrm{dt}\wedge\dAc&=&\int_{\mathrm{S^1}}(\int_{\mathrm{M}}\dAc\wedge\upmu)\mathrm{dt}\\
\nonumber &=&-\frac{i}{2\pi}\int_{\mathrm{S^1}}(\int_{\mathrm{M}}\upnu\wedge\mathrm{d}\Ac)\mathrm{dt}\\
&=&\frac{2\pi}{r}\int_{\mathrm{S^1}}\frac{\mathrm{d}}{\mathrm{d}t}\mathfrak{a}^{\mathcal{F}}(\mathbb{A}_{\underline{\mathbb{C}}},\uppsi_{\underline{\mathbb{C}}})\mathrm{dt}=0.
\end{eqnarray}
Therefore,
\begin{equation}\label{5.53}
\frac{i}{2\pi}\int_{\mathrm{S^1\times M}}\upomega\wedge\FAc=\frac{i}{2\pi}\int_{\mathrm{S^1\times M}}\upomega\wedge\FAct.
\end{equation}
We also note that the left hand side in (\ref{5.53}) is equal to zero since $\Ac$ is a connection on the trivial bundle.\\
\par\emph{Part 2}: Fix $r\geq 1$ large in order to define $\mathfrak{T}_r$ as in Proposition \ref{p3.3}. Let $\mathfrak{T}_r=\{t_{\mathrm{i}}\}_{\mathrm{i}=1,..,\mathrm{N}-r}$. Given $\updelta>0$ very small we shall use $\mathrm{I_{i}}$ to denote the interval $[t_{\mathrm{i}}-\updelta,t_{\mathrm{i}}+\updelta]$ and we shall use $\mathrm{J_{i,i+1}}$ to denote the interval $[t_{\mathrm{i}}+\updelta,t_{\mathrm{i+1}}-\updelta]$. We write the connection $\mathbb{A}_{\mathrm{i,i+1}}$ as $\mathbb{A}_{\mathrm{i,i+1}}=\mathbb{A}_{\mathbb{S}_0}+2\mathrm{A_{i,i+1}}$ where $\mathrm{A_{i,i+1}}$ is viewed as a connection on the bundle $\mathrm{E}$ over $(\mathrm{I_{i}}\cup\bint\cup\mathrm{I_{i+1}})\times\mathrm{M}$. The curvature of $\Aii$ over $\bint\times\mathrm{M}$ is given by
\begin{equation}\label{5.54}
\mathrm{F}_\Aii=\mathrm{F}_\Aiit+\mathrm{dt}\wedge\dot{\mathrm{A}}_{\mathrm{i,i+1}}.
\end{equation}
We now write the integral of $\frac{i}{2\pi}\upomega\wedge(\mathrm{F_\Aii}-\FAct)$ over $\mathrm{\bint\times M}$ as 
\begin{equation} \label{5.28}
\frac{i}{2\pi}\int_{\bint\times\mathrm{M}}\mathrm{dt}\wedge\upnu\wedge(\mathrm{F}_\Aiit-\FAct)+\frac{i}{2\pi}\int_{\bint\times\mathrm{M}}\upmu\wedge\mathrm{dt}\wedge\dot{\mathrm{A}}_{\mathrm{i,i+1}}.
\end{equation}
We will first examine the left most integral in (\ref{5.28}) and then the right most integral. Moreover, in order to consider the left most integral, we fix an integer $n$ to define $\mathrm{J}_{\mathrm{i,i+1};n}$ to be the set of $t\in\bint$ where $\mathcal{E}_{\uptheta}(t)<2^n$. We then consider separately the contribution to the left most integral from $(\bint\setminus\mathrm{J}_{\mathrm{i,i+1};n})\times\mathrm{M}$ and from $\mathrm{J}_{\mathrm{i,i+1};n}\times\mathrm{M}$.\\
\par\emph{Part 3}: Little can be said about the contribution from $(\bint\setminus\mathrm{J}_{\mathrm{i,i+1};n})\times\mathrm{M}$ to the left most integral in (\ref{5.28}) except what is implied by Lemma \ref{l4.1}. In particular, it follows from the latter using (\ref{4.41}) that if $t\in\bint\setminus\mathrm{J}_{\mathrm{i,i+1};n}$, then
\begin{equation}\label{5.101}
\frac{i}{2\pi}\int_{\mathrm{M}_t}\upnu\wedge(\mathrm{F}_\Aiit-\FAct)\geq {c_0}^{-1}\mathcal{E}_\uptheta(t)-c_0
\end{equation}
where $c_0>0$ is independent of $n$, the index $\mathrm{i}$, $t$, and also $r$. Note in particular that (\ref{5.101}) is positive if $2^n>{c_0}^2$.
\par As we show momentarily, there is a positive lower bound for the contribution to the left most integral in (\ref{5.28}) from $\mathrm{J}_{\mathrm{i,i+1};n}\times\mathrm{M}$. To this end, we exhibit constants $c_*>0$ and $r_n>1$ with the former independent of $n$, both independent of $r$ and the index $\mathrm{i}$; and such that
\begin{equation}\label{5.102}
\frac{i}{2\pi}\int_{\mathrm{M}_t}\upnu\wedge(\mathrm{F}_\Aiit-\FAct)\geq c_*
\end{equation}
at each fixed $t\in\mathrm{J}_{\mathrm{i,i+1};n}$ when $r\geq r_n$. What follows is an outline of how this is done. We first appeal to Proposition \ref{p3.7} to find $r_{n}$ such that if $r>r_{n}$, then each point of ${\upalpha_{\mathrm{i,i+1}}}^{-1}(0)$ has distance $c_0r^{-1/2}$ or less from a curve of the vector field that generates the kernel of $\upmu$. We then split the integral in (\ref{5.102}) so as to write it as a sum of two integrals, one whose integration domain consists of points with distance $\mathcal{O}(r^{-1/2})$ or less from the loops in $\mathrm{M}_t$, and the other whose integration domain is complementary part in $\mathrm{M}_t$. We show that the contribution to the former is bounded away from zero by some constant $\mathfrak{L}>0$ which is essentially the length of the shortest closed integral curve of this same vector field. We then show that the contribution from the rest of $\mathrm{M}_t$ is much smaller than this when $r$ is large.\\
\par\emph{Part 4}: Fix  $t\in\mathrm{J}_{\mathrm{i,i+1};n}$. Given $\upepsilon>0$, Proposition \ref{p3.7} finds a constant $r_{n,\upepsilon}$, and if $r>r_{n,\upepsilon}$, a collection $\Theta_t$ of pairs $(\upgamma,\mathrm{m})$ with various properties of which the most salient for the present purposes are that $\upgamma$ is a closed integral curve of the vector field that generates the kernel of $\upmu|_t$ such that $||\upalpha_{\mathrm{i,i+1}}|-|\upmu|^{1/2}|<\upepsilon$ at points with distance $c_\upepsilon r^{-1/2}$ from any loop in $\Theta_t$. Here, $c_{\upepsilon}\geq 1$ depends on $\upepsilon$ but not on $r$, $t$, or the index $\mathrm{i}$. This understood, fix some very small $\upepsilon$ and let $\mathrm{M}_{t,\upepsilon}\subset\mathrm{M}_t$ denote the set of points with distance $2^7c_\upepsilon r^{-1/2}$ or greater from all loops in $\Theta_t$. 
\par To consider the contribution to (\ref{5.102}) from $\mathrm{M}_t\setminus\mathrm{M}_{t,\upepsilon}$, we write the 1-form $\upnu$ as in  (\ref{4.41}). Then, by Lemma \ref{l4.1}, it follows that
\begin{equation}\label{5.103}
\frac{i}{2\pi}\int_{\mathrm{M}_t\setminus\mathrm{M}_{t,\upepsilon}}|\upupsilon\wedge(\mathrm{F}_\Aiit-\FAct)|\leq c_{\upepsilon}r^{-1/2}\mathfrak{L}_t,
\end{equation}
where $\mathfrak{L}_t=\Sigma_{(\upgamma,\mathrm{m})}\mathrm{m}\cdot\mathrm{length}(\upgamma)$.
\par To see about the rest of the $\mathrm{M}_t\setminus\mathrm{M}_{t,\upepsilon}$ contribution, note that Lemma 6.1 in \cite{t1} has a verbatim analogue in the present context. In particular, the latter implies that 
\begin{equation}\label{5.104}
\frac{i}{2\pi}\ast(\ast\upmu\wedge\mathrm{F}_\Aiit)\geq\frac{1}{8\pi}r|\upmu|(|\upmu|-|\upalpha_{\mathrm{i,i+1}}|^2)
\end{equation}
at all points in $\mathrm{M}_t\setminus\mathrm{M}_{t,\upepsilon}$ if $r$ is large. It follows from this, the third item in Proposition \ref{p3.7} and (\ref{5.103}) that 
\begin{equation}\label{5.105}
\frac{i}{2\pi}\int_{\mathrm{M}_t\setminus\mathrm{M}_{t,\upepsilon}}\upnu\wedge(\mathrm{F}_\Aiit-\FAct)\geq c_0\mathfrak{L}_t,
\end{equation}
when $r$ is larger than some constant that depends only on $\upepsilon $ and $n$. Here, $c_0>0$ is independent of $r$, $t$, $n$, $\upepsilon$ and the index $\mathrm{i}$.\\
\par\emph{Part 5}: Turn now to the contribution to (\ref{5.102}) from $\mathrm{M}_{t,\upepsilon}$. By Lemma \ref{l4.2}, no generality is lost by taking $r_{n,\upepsilon}$ so that
\begin{eqnarray}\label{5.106}
\nonumber&&||\upmu|^{1/2}-|\upalpha_{\mathrm{i,i+1}}||<\upepsilon\;\mathrm{and}\;|{\nabla_{\Aii}}^k\upalpha_{\mathrm{i,i+1}}|\leq\upepsilon r^{k/2}\; \mathrm{for}\; k=1,2;\\
&&|{\nabla_{\Aii}}^k\upbeta_{\mathrm{i,i+1}}|\leq\upepsilon r^{(k-1)/2}\; \mathrm{for}\; k=0,1,2
\end{eqnarray}
at all points in $\mathrm{M}_t$ with distance $c_\upepsilon r^{-1/2}$ or more from any loop in $\Theta_t$. Let $\mathrm{M'}$ denote the latter set. Note in this regard that $\mathrm{M}_{t,\upepsilon}$ is the set of points with distance $2^7c_\upepsilon r^{-1/2}$ or more from any loop in $\Theta_t$, so $\mathrm{M}_{t,\upepsilon}\subset\mathrm{M'}$. Meanwhile, we can also assume that (\ref{5.106}) holds at all points in $\mathrm{M}_t$ when $(\Aii,(\upalpha_{\mathrm{i,i+1},\upbeta_{\mathrm{i,i+1}}}))$ is replaced by $(\Ac,(\upalpha_{\underline{\mathbb{C}}},\upbeta_{\underline{\mathbb{C}}}))$. Granted these last observations, we change the gauge for $(\Aii,\uppsi_{\mathrm{i,i+1}})$ on $\mathrm{M'}$ so that $\upalpha_{\mathrm{i,i+1}}=\mathrm{h}\upalpha_{\underline{\mathbb{C}}}$ where $\mathrm{h}$ is a real and positive valued function. Having done so, we write $\Aii$ on $\mathrm{M'}$ as $\Aii=\Ac+(2r)^{1/2}\mathrm{b}$ with $\mathrm{b}$ a smooth imaginary valued 1-form. This understood, then the contribution to (\ref{5.102}) from $\mathrm{M}_{t,\upepsilon}$ is no greater than
\begin{equation}\label{5.107}
c_1\int_{\mathrm{M}_{t,\upepsilon}}|\mathrm{db}|
\end{equation}
where $c_1$ depends only on $\upomega$. Our task now is to show that (\ref{5.107}) is small if $r$ is sufficiently large.
\par To start this task, we note that with our choice of gauge, it follows from (\ref{5.106}) and its $(\Ac,\uppsi_{\underline{\mathbb{C}}})$ analogue that
\begin{equation}\label{5.108}
|\upalpha_{\mathrm{i,i+1}}-\upalpha_{\underline{\mathbb{C}}}|+|\mathrm{b}|\leq c_0\upepsilon
\end{equation}
on $\mathrm{M'}$. Here, $c_0$ is independent of $\upepsilon$ and $r$.
\par Introduce $\mathrm{M''}\subset\mathrm{M'}$ to denote the set of points with distance $2^6c_\upepsilon r^{-1/2}$ or more from any loop in $\Theta_t$. We now see how to find a function $\mathrm{x}:\mathrm{M}\rightarrow\mathbb{R}$ with the following properties: First, $\mathfrak{b}=(\mathrm{b}-i(2r)^{-1/2}\mathrm{dx},e^{i\mathrm{x}}\uppsi-\uppsi_{\underline{\mathbb{C}}},0)$ obeys the equation
\begin{equation}\label{5.109}
\mathcal{L}_{(\Ac,\uppsi_{\underline{\mathbb{C}}})}\mathfrak{b}+r^{1/2}\mathfrak{b}\ast\mathfrak{b}=0
\end{equation}
on $\mathrm{M''}$. Second, $|\mathfrak{b}|\leq\mathrm{z}\upepsilon$ where $\mathrm{z}>0$ is independent of $r$ and $\upepsilon$.
\par To explain our final destination, fix a smooth, non-increasing function $\upchi:[0,\infty)\rightarrow[0,1]$ with value $0$ on $[0,\frac{3}{4}]$ and with value $1$ on $[1,\infty)$. Set ${\upchi_\upepsilon}'$ to denote the function on $\mathrm{M}$ given by
\begin{equation}\label{5.110}
{\upchi_\upepsilon}'=\upchi(dist(\cdot,\cup_{(\upgamma,\mathrm{m})\in\Theta_t}\upgamma)/2^7c_\upepsilon r^{-1/2}).
\end{equation}
Let $\mathfrak{b}'={\upchi_\upepsilon}'\mathfrak{b}$. This function has compact support in $\mathrm{M''}$ and it obeys the equation
\begin{equation}\label{5.111}
\mathcal{L}_{(\Ac,\uppsi_{\underline{\mathbb{C}}})}\mathfrak{b}'+r^{1/2}\mathfrak{b}\ast\mathfrak{b}'=\mathfrak{h},
\end{equation}
where $|\mathfrak{h}|\leq c_0\mathrm{z}|\mathrm{d}{\upchi_\upepsilon}'|\upepsilon $ where $c_0$ is independent of $r$, $t$, $\upepsilon$ and the index $\mathrm{i}$. Note in particular that the $\mathrm{L}^2$-norm of $\mathfrak{h}$ is bounded by $c_1\mathrm{z}\mathfrak{L}_t\upepsilon $ where $c_1$ is also independent of the same parameters. This understood, it follows from (\ref{4.60}) that 
\begin{equation}\label{5.112}
||\mathfrak{b}'||_{\mathbb{H}}\leq c_2\mathrm{z}\upepsilon r^{1/2}||\mathfrak{b}'||_2+c_1\mathrm{z}\upepsilon\mathfrak{L}_t.
\end{equation}
Equation (\ref{5.112}) gives the bound $||\mathfrak{b}'||_{\mathbb{H}}\leq 2c_1\mathrm{z}\upepsilon\mathfrak{L}_t$ when $\upepsilon<\frac{1}{4}(c_2\mathrm{z})^{-1}$. As a final consequence, (\ref{5.107}) is seen to be no greater than $c_3\mathrm{z}\upepsilon\mathfrak{L}_t$ with $c_3$ again independent of $r$, $t$, $\upepsilon$ and the index $\mathrm{i}$.
\par To find the desired function $\mathrm{x}$, introduce again the function $\upchi$, and define $\upchi_{\upepsilon}:\mathrm{M}\rightarrow[0,1]$ by replacing $2^7c_\upepsilon r^{-1/2}$ in (\ref{5.109}) by $2^6c_\upepsilon r^{-1/2}$. Equation (\ref{5.111}) is then satisfied on $\mathrm{M''}$ if $\mathrm{x}$ obeys the equation
\begin{equation}\label{5.113}
\mathrm{d}^{\ast}\mathrm{dx}+2|\upmu|^{1/2}r|\upalpha_{\mathrm{i,i+1}}|\sin{\mathrm{x}}=\upchi_{\upepsilon}\mathrm{d^{\ast}b}.
\end{equation}
\par This equation has the same form as that in (\ref{4.20}). In particular, the arguments in \cite{t2} that find a solution of the equation (2.16) in \cite{t2} can be applied only with minor modifications to find a solution, $\mathrm{x}$, of the equation in (\ref{5.113}) that obeys the bounds in (\ref{5.25}). This being the case, the resulting $\mathfrak{b}=(\mathrm{b}-i(2r)^{-1/2}\mathrm{dx},e^{i\mathrm{x}}\uppsi-\uppsi_{\underline{\mathbb{C}}},0)$ is such that $|\mathfrak{b}|\leq \mathrm{z}\upepsilon$.\\
\par\emph{Part 6}: It follows from what is said in Parts 4 and 5 that there exists $c_*>0$ and $r_n\geq1$ such that if $r\geq r_n$, then (\ref{5.102}) holds. Moreover, $c_*$ is independent of $n$ because it is larger than some fixed fraction of the shortest closed integral curve of any given $t\in\mathrm{S^1}$ version of the kernel of $\upmu$. With (\ref{5.101}), this implies that the left most integral in (\ref{5.28}) obeys 
\begin{equation}\label{5.114}
\frac{i}{2\pi}\int_{\bint\times\mathrm{M}}\mathrm{dt}\wedge\upnu\wedge(\mathrm{F}_\Aiit-\FAct)\geq c_{**}\mathrm{length}(\bint),
\end{equation}
where $c_{**}$ is also independent of $n$ and $r$ which are both very large.
\par To say something about the right most integral in (\ref{5.28}), we write $\Aii=\mathrm{A_E}+\aii$ where $\mathrm{A_E}$ is the $t$-independent connection on $\mathrm{E}$ with harmonic curvature form chosen so that $\mathbb{A_S}=\mathbb{A}_{\mathbb{S}_0}+2\mathrm{A_E}$. We then use the fact that the equations in (\ref{3.1}) are the variational equations of the functional $\mathfrak{a}$ as in (\ref{3.2}) to write

\begin{equation}\label{5.115}
\frac{i}{2\pi}\int_{\mathrm{M}}\upmu\wedge\daii=-\frac{1}{4\pi r}\int_{\mathrm{M}}\aii\wedge\mathrm{d}\aii.
\end{equation}
Here, we use the fact that $\mathcal{D}_{\Aii}\uppsi_{\mathrm{i,i+1}}=0$ to dispense with the derivative of the right most integral in (\ref{3.2}) with respect to $t$. Granted (\ref{5.115}), we identify the right most integral in (\ref{5.28}) with 
\begin{align}\label{5.30}
\nonumber\frac{1}{4\pi r}[-\int_{\mathrm{M}}(\aii\wedge(\mathrm{d}\aii-i\upvarpi_{\mathbb{S}}))|_{t_{\mathrm{i+1}}-\updelta}+\int_{\mathrm{M}}(\aii\wedge(\mathrm{d}\aii-i\upvarpi_{\mathbb{S}}))|_{t_{\mathrm{i}}+\updelta}].\\
\end{align}
Equations (\ref{5.114}) and (\ref{5.30}) summarize what we say for now about (\ref{5.28}).\\
\par\emph{Part 7}: Recall that $\mathrm{I_{i}}=[t_{\mathrm{i}}-\updelta,t_{\mathrm{i}}+\updelta]$. We now review how we define the connection $\mathrm{A_i}$ on $\mathrm{E}$ over $\sint\times\mathrm{M}$. This is done using a `bump' function, $\mathrm{v}:\sint\rightarrow[0,1]$. This function is non-decreasing, it is equal to $0$ near $t_{\mathrm{i}}-\updelta$ and equal to $1$ near $t_{\mathrm{i}}+\updelta$. Meanwhile, we chose gauges for $\mathrm{A_{i-1,i}}$ and $\Aii$ so that there is no spectral flow between the respective $(\mathbb{A}_{\mathrm{i-1,i}},\uppsi_{\mathrm{i-1,i}})$ and $(\mathbb{A}_{\mathrm{i,i+1}},\uppsi_{\mathrm{i,i+1}})$ versions of (\ref{4.18}). Having done so, we write $\mathrm{A_{i-1,i}}=\mathrm{A_E}+\mathrm{a_{i-1,i}}$ and $\mathrm{A_{i,i+1}}=\mathrm{A_E}+\mathrm{a_{i,i+1}}$. We then defined $\mathbb{A}_{\mathrm{i}}=\mathbb{A_S}+2(1-\mathrm{v})\mathrm{a_{i-1,i}}+2\mathrm{v}\mathrm{a_{i,i+1}}$ and we used the latter to define $\Phi$ on $\sint\times\mathrm{M}$ by $\frac{i}{2\pi}(\mathrm{F}_{\mathbb{A}_{\mathrm{i}}}-\mathrm{F}_{\mathbb{A}_{\underline{\mathbb{C}}}})$.
\par In order to say something about 
\begin{equation}\label{5.116}
\int_{\sint\times\mathrm{M}}\upomega\wedge\frac{i}{2\pi}(\mathrm{F}_{\mathbb{A}_{\mathrm{i}}}-\mathrm{F}_{\mathbb{A}_{\underline{\mathbb{C}}}})
\end{equation}
we write $\mathrm{F}_{\mathrm{A_i}}-\FAct$ as 
\begin{eqnarray}\label{5.117}
\nonumber&&\mathrm{v}\;(\mathrm{F}_{\mathrm{A_{i,i+1}}|_t}-\FAct)+(1-\mathrm{v})(\mathrm{F}_{\mathrm{A_{i-1,i}}|_t}-\FAct)\\
&&+\mathrm{dt}\wedge\frac{\partial}{\partial t}(\mathrm{v}\mathrm{a_{i,i+1}})+\mathrm{dt}\wedge\frac{\partial}{\partial t}((1-\mathrm{v})\mathrm{a_{i-1,i}}).
\end{eqnarray}
As we saw in Parts 4 and 5 above, the two left most terms in (\ref{5.117}) give positive contribution to the integral in (\ref{5.116}). The contribution of the two right most terms are
\begin{equation}\label{5.118}
\frac{i}{2\pi}\int_{\sint\times\mathrm{M}}(\mathrm{dt}\wedge\upmu\wedge\frac{\partial}{\partial t}(\mathrm{v}\mathrm{a_{i,i+1}}))+\frac{i}{2\pi}\int_{\sint\times\mathrm{M}}(\mathrm{dt}\wedge\upmu\wedge\frac{\partial}{\partial t}((1-\mathrm{v})\mathrm{a_{i-1,i}})).
\end{equation}
We analyze (\ref{5.118}) using an integration by parts to write it as the sum of
\begin{equation}\label{5.119}
-\frac{i}{2\pi}\int_{\sint\times\mathrm{M}}(\mathrm{dt}\wedge\mathrm{d}\upnu\wedge\mathrm{v}\mathrm{a_{i,i+1}}+(1-\mathrm{v})\mathrm{a_{i-1,i}}),
\end{equation}
and
\begin{equation}\label{5.120}
\frac{i}{2\pi}\int_{\mathrm{M}}(\upmu\wedge\mathrm{a_{i,i+1}})|_{t_{\mathrm{i}}+\updelta}-\frac{i}{2\pi}\int_{\mathrm{M}}(\upmu\wedge\mathrm{a_{i-1,i}})|_{t_{\mathrm{i}}-\updelta}.
\end{equation}
\par Our only remark about the term in (\ref{5.119}) is that it is bounded below by $-\mathcal{K}\updelta$, where $\mathcal{K}$ is a constant that is independent of $\updelta$. This is all we need to know. Meanwhile, we use (\ref{3.2}) to write (\ref{5.120}) as the sum of the two terms:
\begin{equation}\label{5.121}
-\frac{1}{2\pi r}(\mathfrak{a}(\mathfrak{c}_{\uptheta,[t_{\mathrm{i}},t_{\mathrm{i+1}}]})|_{t_{\mathrm{i}}+\updelta}-\mathfrak{a}(\mathfrak{c}_{\uptheta,[t_{\mathrm{i-1}},t_{\mathrm{i}}]})|_{t_{\mathrm{i}}-\updelta})
\end{equation}
and
\begin{equation}\label{5.122}
\frac{1}{4\pi r}[\int_{\mathrm{M}}(\mathrm{a_{i-1,i}}\wedge(\mathrm{d}\mathrm{a_{i-1,i}}-i\upvarpi_{\mathbb{S}}))|_{t_{\mathrm{i}}-\updelta}-\int_{\mathrm{M}}(\aii\wedge(\mathrm{d}\aii-i\upvarpi_{\mathbb{S}}))|_{t_{\mathrm{i}}+\updelta}].
\end{equation}
\par To say something about (\ref{5.121}), recall that we choose the gauges when defining $\mathrm{a_{i-1,i}}$ and $\mathrm{a_{i,i+1}}$ on $\sint\times\mathrm{M}$ so that the spectral flow $\mathcal{F}$ take the same value on  $(\mathbb{A}_{\mathrm{i-1,i}},\uppsi_{\mathrm{i-1,i}})$ and $(\mathbb{A}_{\mathrm{i,i+1}},\uppsi_{\mathrm{i,i+1}})$. As a consequence,
\begin{align}\label{5.123}
\nonumber-\frac{1}{2\pi r}(\mathfrak{a}(\mathfrak{c}_{\uptheta,[t_{\mathrm{i}},t_{\mathrm{i+1}}]})|_{t_{\mathrm{i}}+\updelta}-\mathfrak{a}(\mathfrak{c}_{\uptheta,[t_{\mathrm{i-1}},t_{\mathrm{i}}]})|_{t_{\mathrm{i}}-\updelta})=-\frac{1}{2\pi r}({\mathfrak{a}^{\mathcal{F}}}_{\uptheta}(t_{\mathrm{i}+\updelta})-{\mathfrak{a}^{\mathcal{F}}}_{\uptheta}(t_{\mathrm{i}-\updelta})).\\
\end{align}
Because the function ${\mathfrak{a}^{\mathcal{F}}}_{\uptheta}$ is continuous and piecewise differentiable, what appears on the right hand side of (\ref{5.123}) is bounded below by $-\mathcal{K}\updelta$, with $\mathcal{K}$ again a constant that is independent of $\updelta$.
\par We comment on (\ref{5.122}) in Part 8.\\
\par\emph{Part 8}: The terms in (\ref{5.122}) are fully gauge invariant. This understood, we observe that the term with integral of $\mathrm{a_{i,i+1}}\wedge\mathrm{da_{i,i+1}}$ is identical but for its sign to the right most term in (\ref{5.30}). As the signs are, in fact, opposite, these two terms cancel. Meanwhile, the term with $\mathrm{a_{i-1,i}}\wedge\mathrm{da_{i-1,i}}$ is identical but for the opposite sign, to the left most term in the version of (\ref{5.30}) over the interval $\mathrm{J_{i-1,i;\updelta}}$. Thus, it cancels the latter term. This understood, the sum of the various $\{\bint\}_{\mathrm{i}=1,..,\mathrm{N}_r}$ version of (\ref{5.30}) is exactly minus the sum of the various $\{\sint\}_{\mathrm{i}=1,..,\mathrm{N}_r}$ versions of (\ref{5.122}). Thus, they cancel when we sum up the various contributions to $\int_{\mathrm{S^1\times M}}\upomega\wedge\Phi$. This we now do. In particular, we find from (\ref{5.112}) and from what is said above and in Part 7 that 
\begin{equation}\label{5.124}
\int_{\mathrm{S^1\times M}}\upomega\wedge\Phi\geq 4\pi c_{**}-\mathrm{N}_r\mathcal{K}\updelta
\end{equation}
where $\mathcal{K}$ is a constant that is independent of $\updelta$. Thus, if we take $\updelta>0$ sufficiently small, we see that 
\begin{equation}\label{5.125}
\int_{\mathrm{S^1\times M}}\upomega\wedge\Phi> 0.
\end{equation}
\par\emph{Part 9}: With (\ref{5.125}) understood, our proof of Proposition \ref{p3.9} is complete with a 
suitable idenfication of the class defined by $\Phi$ in $\mathrm{H^2}(\mathrm{M};\mathbb{Z})$.  To this end, remark that it follows from our definition of each $\mathbb{A}_\mathrm{{i,i+1}}$ and each $\mathbb{A}_\mathrm{i}$, that $\Phi$ can be written as $\frac{i}{2\pi}(\mathrm{F}_{\mathbb{A}}-\mathrm{F}_{\mathbb{A}_{\underline{\mathbb{C}}}})$ where $\mathbb{A}$ can be written as $\mathbb{A}_{\mathbb{S}_0}+2\mathrm{A}$ where $\mathrm{A}$ is a connection on a line bundle $\mathrm{E}'$ over $\mathrm{S^1\times M}$ whose first Chern class restricts to each $\mathrm{M}_t$ as that of $\mathrm{E}$. Indeed, $\mathbb{A}$ is defined first on each of $\{\bint\times\mathrm{M}\}_{\mathrm{i}=1,..,\mathrm{N}_r}$ as $\{\mathbb{A}_{\mathrm{i,i+1}}=\mathbb{A}_{\mathbb{S}_0}+2\mathrm{A_{i,i+1}}\}_{\mathrm{i}=1,..,\mathrm{N}_r}$, and then on each of $\{\sint\times\mathrm{M}\}_{\mathrm{i}=1,..,\mathrm{N}_r}$ as $\{\mathbb{A}_{\mathrm{i}}=\mathbb{A}_{\mathbb{S}_0}+2\mathrm{A_E}+2(1-\mathrm{v})\mathrm{a_{i-1,i}}+2\mathrm{v}\mathrm{a_{i,i+1}}\}_{\mathrm{i}=1,..,\mathrm{N}_r}$.  These various connections 
were then glued on the overlaps using maps from $\mathrm{M}$ into $\mathrm{S}^1$.
\par  We write $\mathrm{E}'$ as $\mathrm{E\otimes L}$.  Let $0\in\mathrm{S^1}$ denote any chosen point. Given what was just said, $\mathrm{L}$ over $[0,2\pi)\times\mathrm{M}$ is isomorphic to the trivial bundle. As such, it is obtained from the trivial bundle over $[0,2\pi]\times\mathrm{M}$ by identifying the fiber over $\{2\pi\}\times\mathrm{M}$ with that over $\{0\}\times\mathrm{M}$ using a map $u:\mathrm{M}\rightarrow\mathrm{U}(1)$.
To say more about $\mathrm{L}$, we define for each $t\in\mathrm{S}^1$, a section $\uppsi|_t$ of $\mathbb{S}$ as follows: For any given index $\mathrm{i}\in\{1,..,\mathrm{N}_r\}$, define $\uppsi|_t=\uppsi_{\mathrm{i,i+1}}$ on $\bint\times\mathrm{M}$. We then define $\uppsi$ at $t\in\sint$ to be $\mathrm{v}\uppsi_{\mathrm{i,i+1}}+(1-\mathrm{v})\uppsi_{\mathrm{i-1,i}}$ using the same gauge choices that are used above to define $\mathbb{A}_{\mathrm{i}}$. This done, the pair $(\mathbb{A}=\mathbb{A}_{\mathbb{S}_0}+2\mathrm{A},\uppsi)$ 
defines a pair of connection over $\mathrm{S^1\times M}$ for the line bundle $det(\mathbb{S})\otimes\mathrm{L}^2$ and section of the spinor bundle $\mathbb{S}\otimes\mathrm{L}$. We now trivialize $\mathrm{L}$ over $[0,2\pi)\times\mathrm{M}$ so as to view the restrictions 
to any given $\mathrm{M}_t$ of $(\mathbb{A},\uppsi)$ as defining a smooth map from $[0,2\pi)$ into $\mathcal{C}$. There is then the corresponding 1-parameter family of operators whose $t\in[0,2\pi)$ member is the $(\mathbb{A},\uppsi)|_t$
version of (\ref{4.18}). This family has zero spectral flow.  Indeed, this is the case because $\mathbb{A}$ was defined over $\sint$ by interpolating between $\mathbb{A}_{\mathrm{i-1,i}}$ and $\mathbb{A}_{\mathrm{i,i+1}}$ in gauges where there is zero spectral flow between the respective $(\mathrm{A}_{\mathrm{i-1,i}},\uppsi_{\mathrm{i-1,i}})$ and $(\mathrm{A}_{\mathrm{i-1,i}},\uppsi_{\mathrm{i-1,i}})$ versions of (\ref{4.18}).
\par Because $(\mathbb{A},\uppsi)|_{2\pi}=(\mathbb{A}|_0-2u^{-1}\mathrm{d}u,u\uppsi|_0)$ and there is no spectral flow 
between the respective $(\mathbb{A},\uppsi)|_0$ and $(\mathbb{A},\uppsi)|_{2\pi}$ versions of (\ref{4.18}), it follows from \cite{aps} that the cup product of $c_1(\mathrm{L})$ with 
$c_1(det(\mathbb{S}))$ is zero.  
 \par Keeping this last point in mind, and given that $\mathrm{L}$ restricts as the trivial bundle to 
each $\mathrm{M}_t$, we use the K\"unneth formula to see that the cup product of $c_1(\mathrm{L})$ with the class defined by $\upomega$ is the same as that between $c_1(\mathrm{L})$ and the class defined by $\upmu|_{0}$.  By 
assumption, the latter class is proportional to $c_1(det(\mathbb{S}))$. Thus, $c_1(\mathrm{L})$ has zero cup product with $[\upomega]$. \QED

\end{document}